\documentclass[11pt,reqno]{amsart}
\usepackage[utf8]{inputenc}
\usepackage[letterpaper, portrait, margin=1in]{geometry}
\usepackage[mathscr]{euscript}
\usepackage{graphicx}
\usepackage{accents}
\usepackage{amsmath}
\usepackage{amssymb}
\usepackage{amsthm}
\usepackage{enumerate}
\usepackage{enumitem}
\usepackage{tikz}
\usepackage{mathtools}
\usepackage{esint}
\usepackage{multirow}
\usepackage{upgreek}
\usepackage{esvect}
\usepackage{setspace}
\usepackage{bbm}
\usepackage{fancyhdr}

\fancypagestyle{firstpage}{
  \fancyhf{}
  \fancyfoot[L]{%
    \begin{minipage}{\linewidth}
      \rule{1.5in}{0.4pt}\\[4pt]
      \scriptsize
        Department of Mathematics, Louisiana State University, Baton Rouge\\
        Email: \texttt{cbunti5@lsu.edu}
    \end{minipage}%
  }

}

\author{Christopher Bunting}
\title[Ergodicity for NSNE with Shear-Thickening Viscosity]{Ergodicity for the Two-Dimensional Stochastic Navier-Stokes Equations with Shear-Thickening Viscosity}
\date{}

\newcommand{\Lim}[1]{\raisebox{0.5ex}{\scalebox{0.8}{$\displaystyle \lim_{#1}\;$}}}

\newcommand{\RN}[1]{\textup{\uppercase\expandafter{\romannumeral#1}}}

\usepackage{cleveref}

\theoremstyle{plain}
\newtheorem{lemma}{Lemma}[section]

\theoremstyle{plain}
\newtheorem{corollary}{Corollary}[section]

\theoremstyle{definition}

\theoremstyle{plain}
\newtheorem{theorem}{Theorem}[section]

\theoremstyle{plain}
\newtheorem*{nolemma}{Lemma}

\theoremstyle{plain}
\newtheorem{proposition}{Proposition}[section]

\theoremstyle{definition}
\newtheorem{defn}{Definition}[section]

\theoremstyle{plain}
\newtheorem{assumption}{Assumption}[section]

\theoremstyle{remark}
\newtheorem{remark}{Remark}

\makeatletter
\newcommand{\reqnomode}{\tagsleft@false}
\makeatother

\numberwithin{equation}{section}

\makeatletter
\renewcommand{\subsection}{%
  \@startsection{subsection}{2}{\z@}%
    {-2.5ex\@plus-1ex \@minus-.2ex}%
    {1ex}%
    {\normalfont\bfseries}}
\makeatother

\makeatletter
\renewcommand{\@seccntformat}[1]{%
  \text{\csname the#1\endcsname.}\quad}
\makeatother

\usepackage{etoolbox}

\makeatletter
\patchcmd{\subsubsection}
  {\normalfont}
  {\normalfont\bfseries}
  {}{}
\makeatother

\begin{document}

\begin{abstract}
    The two-dimensional stochastic Navier-Stokes equations with a shear-thickening viscosity driven by multiplicative noise are studied. After the problem is formulated as an abstract stochastic evolution equation on suitable function spaces, the existence and uniqueness of strong solutions are established through local monotonicity methods. Moment estimates are then obtained and used to prove the existence of an ergodic invariant probability measure via the Krylov-Bogolyubov technique. The uniqueness of the ergodic invariant measure is subsequently established by constructing a generalized coupling through a controlled auxiliary equation. Finally, the results are extended to other types of viscosities as well as to different stochastic hydrodynamical models.
\end{abstract}

\maketitle

\thispagestyle{firstpage}



\section{Introduction}

Let $G \subset \mathbb{R}^2$ be a bounded domain with a smooth boundary $\partial G$. Let $u$ denote the velocity, $\rho$, the pressure field, and $f$, the deterministic external forcing. Consider the 2D Navier-Stokes equation for incompressible fluids with nonlinear viscosity and no-slip boundary condition perturbed by a stochastic external force, i.e.
\begin{align}
    \frac{\partial u}{\partial t} - \nu_0 \Delta u - \nu_1\nabla \cdot (|\nabla u|^{q-2} \nabla u) + (u \cdot \nabla) u + \nabla \rho &= f(t) + \sigma(u) \frac{dW(t)}{dt} \notag\\
    \nabla\cdot u &= 0 \label{eq:snse}
\end{align}
with $ u(x,t) = 0$ for $ x\in \partial G$ and $u(x,0) = u_0(x)$ for all $x\in G$.\\

Here, $\sigma(u)$ denotes a suitable multiplicative noise coefficient, $\nu_0 > 0$, the kinematic viscosity, $\nu_1 >0$, the strength of the nonlinear viscosity operator, and $q > 2$ is related to the shear-thickening behavior. When $\nu_1 = 0$ or $q =2$, the system\;\eqref{eq:snse} reduces to the well-studied stochastic Navier-Stokes equation perturbed by a stochastic force.\\

The objective of the present article is to establish the existence of an ergodic  measure for the equation\;\eqref{eq:snse} in the presence of shear-thickening viscosity and multiplicative noise.  The idea of introducing a Wiener process in the Navier–Stokes equation in order to study invariant measures was an idea originally proposed by Kolmogorov (see\;\cite{VF}) so that the problem of stationarity of distributions for the solution can be studied, and such ideas can also be found in the classical book by Landau and Lifshitz\;\cite{LL}. Early mathematical work on invariant measures for stochastic fluid equations can be found in Albeverio and Cruzeiro\;\cite{AC}, who considered invariant Gibbs measures for two-dimensional Euler flows and the corresponding stochastic Navier-Stokes dynamics. Since then, a substantial literature has developed concerning invariant measures and ergodicity for stochastic Navier-Stokes equations under a variety of stochastic perturbations. For the two-dimensional equations, invariant measures and ergodicity have been studied under Wiener and random forcing in\;\cite{Fl,FG,FM,Ma,MH,sundar}, including the case of degenerate stochastic forcing\;\cite{MH}. Other types of noise have also been considered, including space-time white noise\;\cite{BF}, fractional Brownian noise\;\cite{FSV}, and multiplicative Wiener noise\;\cite{ferr}. In the three-dimensional setting, ergodicity has been studied for equations driven by Gaussian Wiener noise\;\cite{DD} and, more recently, by L\'{e}vy noise\;\cite{mss}. These works demonstrate the broad range of stochastic perturbations for which invariant measures and their long-time statistical properties have been investigated.\\

The equation considered in the present paper differs from the classical stochastic Navier-Stokes equation through the presence of the nonlinear viscosity
\[-\nu_1\nabla\cdot\left(|\nabla u|^{q-2}\nabla u\right), \qquad q>2.\]
Such nonlinear viscosity models arise naturally in the mathematical description of non-Newtonian fluids. In particular, the case $q>2$ corresponds to shear-thickening behavior, for which the effective viscosity increases with the magnitude of the shear. Nonlinear viscosity models of this type were introduced in the deterministic setting by Ladyzhenskaya\;\cite{oa1,oa2}, and the mathematical theory of nonlinear parabolic equations and power-law fluids has been developed extensively; see, for example,\;\cite{barbu,dib,jll,jll2}. The stochastic counterpart has also been studied. In particular, Breit\;\cite{Breit} established the existence of martingale weak solutions for stochastic power-law fluids driven by Brownian noise. These results provide an important foundation for the analysis of stochastic fluid equations with nonlinear constitutive laws.\\

The novelty in this work arises from the presence of this nonlinear viscosity term in the Navier-Stokes system along with an external body force in the form of a continuous multiplicative noise. The solvability of the nonlinear stochastic Navier-Stokes equation (NSNE, for short), the existence of invariant measures for  NSNE and its uniqueness are established in this paper along with an extension of our results to the three-dimensional convolution-regularized NSNE.\ After casting the NSNE as an abstract evolution equation on a suitable Hilbert space, the existence and uniqueness of the solution is shown by the method of local monotonicity\;\cite{men}. The existence of invariant measures is proved by the Krylov-Bogolyubov argument. Next, by introducing an oppositional control term in the NSNE\;\cite{ghmr}, a related equation called the stopped nudged equation is shown to be equivalent in law to the NSNE. The path-equivalence of solutions of the original equation to those of the transformed equation allows one to conclude ergodicity for the NSNE  when the noise coefficient is chosen in a class of integrands (see\;\cite{ferr, ghmr}).\\

It is worthwhile to mention that our results apply to a class of stochastic partial differential equations that includes the stochastic MHD system with a shear-thickening term, and analogous stochastic shell models in turbulence theory such as the GOY model. The results can be modified for stochastic models with other nonlinear viscosities such as a shear-dependent viscosity. The study of ergodic behavior for solutions of the stochastic Navier-Stokes equations with general viscosities and suitable multiplicative noise driven by Poisson random measures is under preparation.\\

The paper is organized as follows. In section 2, the functional analytic setup and background results are given. A priori estimates and solvability of NSNE is presented in section 3. Existence of invariant measures for the system is given in section 4. Ergodic behavior of the stochastic system is established in section 5. Extensions of our results to other models and nonlinear viscosities are outlined and discussed in section 6.

\section{Functional Setting and Assumptions}
Equation\;\eqref{eq:snse} can be written in an abstract evolution form by introducing the following function spaces and operators. Let $\mathcal{V}$ be the space $C^\infty$ functions satisfying $\nabla \cdot u = 0$ and the boundary condition. Define the spaces $H$ as the completion of $\mathcal{V}$ in the $L^2(G)$ norm, $V$ as the completion of $\mathcal{V}$ in the $W^{1,2}_0(G)$ norm, and $V_{1,q}$ as the completion of $\mathcal{V}$ in the $W^{1,q}_0(G)$ norm. The following characterization of $H$, $V$, and $V_{1,q}$ is known\;\cite{tm}.
\begin{align*}
    H &= \{u\in L^2(G); \nabla \cdot u = 0, u \cdot \textbf{n}|_{\partial G} = 0\}\\
    V &= \{u \in W^{1,2}_0(G); \nabla \cdot u =0\}\\
    V_{1,q}&= \{u \in W^{1,q}_0(G); \nabla \cdot u = 0\}.
\end{align*}
where \textbf{n} is the outward normal and $u\cdot \textbf{n}|_{\partial G}$ is in the sense of trace operators\;\cite{cf}. Let $H'$, $V'$, and $V_{1,q}'$ denote the duals of $H$, $V$, and $V_{1,q}$ respectively. Since $G$ is a bounded domain, the norms of these spaces satisfy the following equivalences
\begin{align*}
    ||u||^2_H &= \int_G |u|^2 dx\\
    ||u||^2_V &= \int_G |\nabla u|^2 dx\\
    ||u||_{V_{1,q}}^q &= ||\nabla u||_{L^q}^q = \int_G |\nabla u|^q dx. 
\end{align*}
For simplicity of notation, we will denote the norm of $H$ by $|u|$ and the norm of $V$ by $||u||$. These spaces form the extended Gelfand triple
\begin{equation}\label{eq:embedding}
    V_{1,q} \hookrightarrow V \hookrightarrow H = H' \hookrightarrow V'\hookrightarrow V_{1,q}'
\end{equation}
where the embeddings are compact, dense, and continuous.\\

We define the operator $A : V \rightarrow V'$ by $Au = -\mathbf{P}_H \Delta u$ for $u\in D(A)$ where $\mathbf{P}_H$ denotes the orthogonal projection into $H$ and the domain $D(A) = \{u \in H : -\mathbf{P}_H\Delta u \in H\}$. By the Cattabriga-Solonnikov regularity theorem (see Chapter 3 of\;\cite{oa2}), $D(A) = W^{2,2}(G) \cap V$. The operator $A$ is the Stokes operator and is a positive, self-adjoint operator. We will denote the eigenfunctions of $A$ as ${\{e_k\}}_{k=1}^{\infty}$ and their corresponding eigenvalues as ${\{\lambda_k\}}_{k=1}^{\infty}$ which satisfy $0< \lambda_1 < \lambda_2 \leq \lambda_3 \cdots \leq \lambda_k \leq \cdots$ and $\lambda_k \rightarrow \infty$ as $k \rightarrow \infty$. Using integration by parts, we say
\begin{equation}\label{eq:aw}
  \langle Au , u \rangle = \langle A^{1/2}u, A^{1/2} u\rangle = |A^{1/2}u|^2 = ||u||^2.  
\end{equation}
where the last equality is due to $G$ being a bounded domain in $\mathbb{R}^2$. Although this shows that the $V$-norm can be expressed in terms of the $H$-norm, we will continue using $||\cdot ||$ to denote the norm of $V$.\\

{Define $\mathcal{A}_1(u)$, for each $u\in V_{1,q}$, such that 
\begin{equation}\label{eq:weaknl}
    \langle \mathcal{A}_1(u), v\rangle = \int_G |\nabla u|^{q-2} \nabla u\; \nabla v\; dx \quad \forall v\in V_{1,q}.
\end{equation}
It is clear that $\langle \mathcal{A}_1(u), v\rangle$ is linear in $v$. Since $u,v \in V_{1,q}$, we have that $\nabla u, \nabla v \in L^q$ and $|\nabla u|^{q-2} \nabla u \in L^{\frac{q}{q-1}}$. It follows that, by the H\"{o}lder inequality,
\begin{align*}
    |\langle \mathcal{A}_1(u), v \rangle| &= \bigg|\int_G |\nabla u|^{q-2} \nabla u\; \nabla v\; dx\bigg|\\
    &\leq ||\nabla u||_{L^q}^{q-1}\; ||\nabla v||_{L^q} < \infty.
\end{align*}
So, the integral $\langle \mathcal{A}_1(u), v\rangle$ is well-defined. Additionally, we observe that, since $||\nabla v||_{L^q} = ||v||_{V_{1,q}}$, $\mathcal{A}_1(u)$ is bounded, i.e.
\begin{equation}\label{eq:nlbdd}
  ||\mathcal{A}_1(u)||_{V'_{1,q}} \leq ||u||_{V_{1,q}}^{q-1}.  
\end{equation}
Thus, the existence of such an operator $\mathcal{A}_1(u)$ satisfying\;\eqref{eq:weaknl} is ensured. Now, we say that the nonlinear operator $\mathcal{A}_1 : V_{1,q}\rightarrow V_{1,q}'$ is defined such that $\mathcal{A}_1(u) = -\mathbf{P}_H \nabla \cdot (|\nabla u|^{q-2}\nabla u)$ where this equality holds in $V'_{1,q}$, i.e.\; in the weak (variational) sense. Similar to the Stokes operator, we say
\[\langle \mathcal{A}_1(u) ,u \rangle = ||\nabla u||_{L^q}^q = ||u||_{V_{1,q}}^q.\] 

Next, we introduce an important lower bound and some properties of $\mathcal{A}_1$.
\begin{proposition}\label{prop:lowerbound}
    Let $u\in V_{1,q}$ and $p\geq 2$. Then 
    \[|u|^{p-2}||u||_{V_{1,q}}^q \geq \tilde{C}(|u|^p - 1)\]
    where $\tilde{C} = \lambda_1^{q/2} |G|^{1-q/2}$ where $|G|$ is the size of the bounded domain.
\end{proposition}
\begin{proof}
    First, note that 
    \[|u|^{p-2} ||u||_{V_{1,q}}^q = |u|^{p-2} ||\nabla u||_{L^q}^q.\]
    By the Poincar\'{e} inequality, we have that 
    \[|u|^{p-2} ||u||_{V_{1,q}}^q \geq C |u|^{p-2} ||u||_{L^q}^q.\]
    Since $q> 2$ and $G$ is bounded, $L^q(G) \subseteq L^2(G)$ and hence 
    \[|u|^{p-2} ||u||_{V_{1,q}}^q \geq C |u|^{p-2} ||u||_{L^q}^q \geq \tilde{C} |u|^{p-2} |u|^q = \tilde{C} |u|^{p+q-2}.\]
    Here $\tilde{C} = \lambda_1^{q/2} |G|^{1-q/2}$ where $|G|$ is the size of the bounded domain. To arrive at the desired lower bound, we use the following lemma.
    \begin{nolemma}
        Let $\alpha \geq 0$ and $p,q \geq 2$. Then 
        \[\alpha^{p+q-2} \geq \alpha^p -1.\]
    \end{nolemma}
    \begin{proof}
        We separate this proof into two cases: $0\leq \alpha \leq 1$ and $\alpha \geq 1$. Suppose $0\leq \alpha\leq 1$. It is clear that $0\leq \alpha^{p+q-2},\alpha^{p} \leq 1$. Hence 
        \[\alpha^{p} - 1 \leq 0 \leq \alpha^{p+q-2}.\]
        The first case is complete. Now, suppose $\alpha \geq 1$. It is clear that $\alpha^p \geq 1$. Since $q\geq 2$, we know that $p+q-2 \geq p$. It follows that 
        \[\alpha^{p+q-2} \geq \alpha^p \geq \alpha^p -1.\]
        The second case is complete. Thus, the inequality holds for any nonnegative $\alpha$. 
    \end{proof}
    \noindent
    Using this lemma, we have 
    \[|u|^{p-2} ||u||_{V_{1,q}}^q \geq \tilde{C} |u|^{p+q-2} \geq \tilde{C}(|u|^p -1).\]
\end{proof}
\noindent
The following are properties of the nonlinear viscosity operator.
\begin{proposition}[\cite{s}, Theorem 1]
    The nonlinear viscous operator $\mathcal{A}_1$ is strictly monotone, i.e.
    \begin{equation}\label{eq:nl}
        \langle \mathcal{A}_1(u) - \mathcal{A}_1(v) , u-v \rangle \geq \gamma||u-v||_{V_{1,q}}^q.  
    \end{equation}
    Here $\gamma$ depends only on the dimension and $q$.
\end{proposition}}
{
\begin{proof}
    Note that for $u,v$
    \[\langle \mathcal{A}_1(u) - \mathcal{A}_1(v), u-v \rangle = \int_G (|\nabla u|^{q-2}\nabla u - |\nabla v|^{q-2}\nabla v) \cdot (\nabla u - \nabla v)dx.\]
    This integrand will be estimated from below using the following lemma.
    \begin{nolemma}[\cite{dib}, Section 1.4]
        Let $q\geq 2$. Then, $\forall x,y\in \mathbb{R}^n$, we have
        \[(|x|^{q-2}x-|y|^{q-2}y) \cdot (x-y) \geq \gamma(q,n) |x-y|^q.\]
    \end{nolemma}
\noindent
Using this lemma, we get 
    \[\langle \mathcal{A}_1(u) - \mathcal{A}_1(v), u-v \rangle \geq \gamma\int_G |\nabla u - \nabla v|^q dx = \gamma||\nabla u - \nabla v||^q_{L^q} = \gamma||u-v||_{V_{1,q}}^q\]
    as desired.
\end{proof}}
{
\begin{proposition}[\cite{s}, Theorem 2]
    The operator $\mathcal{A}_1(\cdot)$ is demicontinuous, i.e.\; let $u_n \rightarrow u$ in $V_{1,q}$, then $\mathcal{A}_1(u_n) \rightarrow \mathcal{A}_1(u)$ in the weak-star topology of $V'_{1,q}$.
\end{proposition}}
{
\begin{proof}
    Let $u_n \rightarrow u$ strongly in $V_{1,q}$. Then, by\;\eqref{eq:nlbdd}, we can take a subsequence (again denoting) $u_n \in V_{1,q}$ such that $\mathcal{A}_1(u_n) \rightarrow \kappa$ weakly in $V'_{1,q}$. We now show that $\kappa = \mathcal{A}_1(u)$. In fact, $u_n \rightarrow u$ strongly in $W^{1,q}(G)$ implies $\nabla u_n \rightarrow \nabla u$ strongly in $L^q(G)$ and hence $\nabla u_n(x) \rightarrow \nabla u(x)$, $x$ a.e.\; in $G$. Hence, $|\nabla u_n(x)|^{q-2}\nabla u_n(x) \rightarrow |\nabla u(x)|^{q-2}\nabla u(x)$, $x$ a.e.\; in $G$. By Lemma 1.3 in Section 1.4 of\;\cite{jll2}, we conclude that, for all $v\in V_{1,q}$,
    \[\langle \mathcal{A}_1(u_n), v\rangle = \int_G |\nabla u_n|^{q-2}\nabla u_n \; \nabla v\; dx \rightarrow \int_G |\nabla u|^{q-2}\nabla u \; \nabla v\; dx.\]
\end{proof}}
{
\begin{corollary}[\cite{s}, Corollary 1]
    The operator $\mathcal{A}_1(\cdot)$ is maximal monotone.
\end{corollary}}
{Maximal monotonicity is a result of demicontinuity and monotonicity, see\;\cite{barbu}. The following corollary is proved in Theorem 2.1 of Chapter 2 Section 2.1 of\;\cite{jll2}. 
\begin{corollary}[\cite{s}, Corollary 2]
    The operator $\mathcal{A}_1 : V_{1,q} \rightarrow V'_{1,q}$ is surjective.
\end{corollary}}
{
\begin{corollary}[\cite{s}, Corollary 3]\label{corollary:hemi}
    The operator $\mathcal{A}_1(\cdot)$ is strongly continuous on finite dimensional subspaces. That is, for all $V_N \subset V_{1,q}$ where $V_N$ is finite-dimensional, the map $u\rightarrow \mathcal{A}_1(u)$ is continuous from $V_N \rightarrow V'_{1,q}$.
\end{corollary}}
\noindent
\Cref{corollary:hemi} follows from Lemma 1.1 in Section 1 of Chapter II of\;\cite{barbu}.\\

Now, we define the trilinear operator $b : V \times V \times V \rightarrow \mathbb{R}$ by 
\[b(u,v,w) = \sum_{i,j=1}^2 \int_G u_i \frac{\partial v_j}{\partial x_i} w_j \, dx.\]
Using this definition, we define the continuous bilinear operator $B: V \times V \rightarrow V'$ such that 
\[\langle B(u,v), w \rangle = b(u,v,w)\]
for all $u, v, w \in V$. We will use $B(u)$ to denote $B(u,u)$. Note that, since $u\in V$,
\begin{equation}\label{eq:tl}
    b(u,v,w) = - b(u,w,v)
\end{equation}
and it follows that $b(u,v,v) = 0$. We also have the following estimates on the bilinear term.
\begin{proposition}
    If $u,v,w \in V$, the bilinear operator $B$ satisfies 
    \begin{equation}\label{eq:bl}
        |\langle B(u,v) , w\rangle|^2 \leq \eta^2 |u|\, ||u||\, ||v||^2\, |w|\, ||w||.
    \end{equation}
\end{proposition}

\begin{proof}
    Let $u,v,w \in V$. Consider $|\langle B(u,v),w\rangle|$. Applying the H\"{o}lder inequality, we have 
    \[|\langle B(u,v),w \rangle| \leq ||u||_{L^4}\, ||v||\, ||w||_{L^4}.\]
    The $L^4$ norms can be estimated by the Ladyzhenskaya inequality, i.e. $||u||^2_{L^4} \leq \eta|u|\, ||u||$, to get 
    \[|\langle B(u,v),w \rangle| \leq \eta|u|^{1/2}\, ||u||^{1/2}\, ||v||\, |w|^{1/2}\, ||w||^{1/2}.\]
    Squaring this estimate yields the desired result.
\end{proof}

We now consider $B(\cdot,\cdot)$ restricted to $V_{1,q}\times V_{1,q}$. That is, $B: V_{1,q} \times V_{1,q} \rightarrow V_{1,q}'$. Under this restriction, we get the following estimate.
\begin{proposition}
    If $u,v,w \in V_{1,q}$, the bilinear operator $B$ satisfies 
    \begin{equation}\label{eq:bl2}
        |\langle B(u,v) , w\rangle|^2 \leq C|u|\, ||u||\, ||v||_{V_{1,q}}^2\, |w|\, ||w||.
    \end{equation}
\end{proposition}

\begin{proof}
    Let $u,v,w\in V_{1,q}$. Applying the H\"{o}lder inequality to $|\langle B(u,v),w\rangle|$, we have 
    \[|\langle B(u,v),w\rangle| \leq ||u||_{L^4}\, ||v||_{V_{1,q}}\, ||w||_{L^r}.\]
    Here $r$ satisfies \(\frac{1}{4} + \frac{1}{q} + \frac{1}{r} = 1\). Since $q>2$, we have that $r < 4$. Using the Ladyzhenskaya inequality on $||u||_{L^4}$, we say 
    \[|\langle B(u,v),w\rangle| \leq \eta^{1/2}|u|^{1/2}\,||u||^{1/2}\, ||v||_{V_{1,q}}\, ||w||_{L^r}.\]
    Since \(r<4\), we interpolate between $L^2$ and $L^4$ so that, for a suitable \(\theta\in (0,1)\), 
    \[|\langle B(u,v),w\rangle| \leq \eta^{1/2}|u|^{1/2}\,||u||^{1/2}\, ||v||_{V_{1,q}}\, |w|^{\theta}\, ||w||_{L^4}^{1-\theta}.\]
    Then, using the fact that \(||w||_{L^4} \leq C||w||\) yields 
    \[|\langle B(u,v),w\rangle| \leq C|u|^{1/2}\,||u||^{1/2}\, ||v||_{V_{1,q}}\, |w|^{\theta}\, ||w||^{1-\theta}\]
    for a suitable constant $C$. Finally, for a possibly different constant $C$, we may write 
    \[|\langle B(u,v),w\rangle|\leq C|u|^{1/2}\,||u||^{1/2}\, ||v||_{V_{1,q}}\, |w|^{1/2}\,||w||^{1/2}.\]
    Squaring this inequality yields the desired result.
\end{proof}
From estimate\;\eqref{eq:bl2}, we conclude that, for $u\in V_{1,q}$, $B(u) \in V_{1,q}'$ and 
\[||B(u)||^2_{V_{1,q}'} \leq C|u|\, ||u||.\]
The deterministic forcing $f$ is assumed to be $V'$-valued and belong to $L^2(0,T; V')$, i.e.
\[||f||_{L^2(0,T; V')}^2 = \int_0^T ||f||_{V'}^2 < \infty.\]
Let $Q$ be a positive, symmetric, trace-class operator on $H$.
\begin{defn}
    A stochastic process $\{W(t) : 0\leq t\leq T\}$ is said to be an $H$-valued $\{\mathcal{F}_t\}$-adapted Wiener process with nuclear covariance operator $Q$ if 
    \begin{enumerate}[label= (\arabic*)]
        \item for each nonzero $h\in H$, $|Q^{1/2}h|^{-1}(W(t),h)$ is a standard one-dimensional Wiener process, and
        \item for any $h\in H$, $(W(t), h)$ is a martingale adapted to $\{\mathcal{F}_t\}$.  
    \end{enumerate}
\end{defn}
\smallskip
Let $U = Q^{1/2} H$ with inner product
\[{(u,v)}_{U} = (Q^{-1/2} u, Q^{-1/2} v)\]
for all $u,v\in U$. Then $U$ is a Hilbert space. Let $||\cdot ||_{U}$ denote the norm in $U$. The embedding $U$ into $H$ is Hilbert-Schmidt since $Q$ is a trace class operator. Let $L_Q(U,H)$ denote the space of linear operators $S$ such that $SQ^{1/2}$ is a Hilbert-Schmidt operator from $H$ to $H$. Define the norm of $L_Q$ by $||S||_{L_Q} = \text{tr}[SQS^*]$.
\begin{assumption}\label{ass:cv}
Assume that the noise coefficient $\sigma : H \rightarrow L_Q(U, H)$ satisfies the following hypotheses.
\begin{enumerate}[label = \textbf{(M\arabic*)}]
    \item\label{M1} The function $\sigma \in C(H; L_Q(U, H))$.
    \item\label{M2} There exists a constant $K>0$ such that for all $u\in H$, $||\sigma(u)||_{L_Q}^2 \leq K(1 + |u|^2)$.
    \item\label{M3} There exists a constant $L>0$ such that for all $u,v\in H$, $||\sigma(u) -\sigma(v)||^2_{L_Q} \leq L|u-v|^2$
\end{enumerate}
\end{assumption}

Applying the Leray projection to the nonlinear stochastic Navier-Stokes equations, we write the system\;\eqref{eq:snse}, in differential form, as 
\begin{equation}\label{eq:diffform}
    \frac{du}{dt} + \nu_0 Au + \nu_1\mathcal{A}_1(u) + B(u) = f(t) +\sigma(u) \frac{dW(t)}{dt}
\end{equation}
in the sense that $u$ solves the following integral equation
\begin{equation}\label{eq:af}
    u(t) + \int_0^t (\nu_0 Au + \nu_1\mathcal{A}_1(u) + B(u))ds = u_0 + \int_0^t f(s) ds + \int_0^t \sigma(u) dW(s).
\end{equation} 
With the estimates above, we provide the following result concerning local monotonicity.

\begin{proposition}\label{prop:localmono}
    For a given $r>0$, let $B_r$ denote the $L^4(G)$ ball in $V_{1,q}$, i.e. $B_r = \{v\in V_{1,q} : ||v||_{L^4(G)} \leq r\}$. Define the nonlinear operator $F$ by $F(u) = -\nu_0Au -\nu_1\mathcal{A}_1(u) - B(u)$. Then, for $2L<\nu_0$, the nonlinear operator $F$ is monotone in $B_r$ i.e., for any $u\in V_{1,q}$ and $v\in B_r$, if $w$ denotes $u-v$, then 
    \begin{equation}
        \langle F(u) - F(v), w \rangle - \frac{r^4}{\nu_0^3}|w|^2 + ||\sigma(u) - \sigma(v)||^2_{L_Q} \leq 0.
    \end{equation}
\end{proposition}

\begin{proof}
    Using\;\eqref{eq:tl} and the bilinearity of the operator $B$, it follows that
    \[\langle B(u) , w\rangle = \langle B(u), u\rangle - \langle B(u), v\rangle + \langle B(u,v),v \rangle = - \langle B(u,w),v\rangle.\]
    Similarly, $\langle B(v) , w\rangle = - \langle B(v,w),v\rangle$. Using these identities, we obtain
    \[\langle B(u) - B(v) , w\rangle = -\langle B(w), v\rangle.\]
    Applying the H\"{o}lder and Ladyzhenskaya inequalities, we arrive at 
    \begin{equation}\label{eq:bw}
    \begin{aligned}[b]
        |\langle B(w),v \rangle| &\leq ||w||_{L^4} ||w|| \, ||v||_{L^4}\\
        &\leq |w|^{1/2} ||w||^{1/2} ||w|| \, ||v||_{L^4}\\
        &= |w|^{1/2}||w||^{3/2} \, ||v||_{L^4}\\
        &\leq \frac{\nu_0}{2} ||w||^2 + \frac{27}{32\nu_0^3} |w|^2 ||v||^4_{L^4}
    \end{aligned}
    \end{equation}
    where the last inequality is a consequence of the basic inequality due to Young. Combining\;\eqref{eq:aw},\;\eqref{eq:nl}, and\;\eqref{eq:bw} yields
    \begin{align*}
        \langle F(u) - F(v), w \rangle & = \langle -\nu_0 Aw , w \rangle + \langle -\nu_1 \mathcal{A}_{1} (u)  + \nu_1 \mathcal{A}_{1}(v), w \rangle - \langle B(u) - B(v), w \rangle\\
        &\leq -\nu_0 ||w||^2 -\nu_1 \gamma ||\nabla w||^q_{L^q} + \frac{\nu_0}{2}||w||^2 + \frac{27}{32\nu_0^3}|w|^2 ||v||^4_{L^4}\\
        &= \frac{-\nu_0}{2}||w||^2 - \nu_1 \gamma ||w||^q_{V_{1,q}} + \frac{27}{32\nu_0^3}|w|^2 ||v||^4_{L^4}\\
        &\leq \frac{-\nu_0}{2}||w||^2 - \nu_1 \gamma ||w||^q_{V_{1,q}} + \frac{r^4}{\nu_0^3}|w|^2.
    \end{align*}
    Note that, by Assumption\;\ref{M3},
    \begin{align*}
        \frac{1}{L}||\sigma(u) - \sigma(v)||^2_{L_Q} \leq |w|^2 \leq ||w||^2. 
    \end{align*}
    From this, we say
    \begin{align*}
      \langle F(u) - F(v), w \rangle - \frac{r^4}{\nu_0^3}|w|^2 + \frac{\nu_0}{2L}||\sigma(u) - \sigma(v)||^2_{L_Q} \leq 0  
    \end{align*}
    Additionally, if $2L<\nu_0$, then
    \begin{align*}
      \langle F(u) - F(v), w \rangle - \frac{r^4}{\nu_0^3}|w|^2 + ||\sigma(u) - \sigma(v)||^2_{L_Q} \leq 0  
    \end{align*}
\end{proof}

Let $H_n := \text{span}\{e_1, e_2, \dots , e_n\}$ where $\{e_j\}$ is a fixed orthonormal basis of eigenfunctions of $A$ in $H$. Let $P_n$ denote the orthogonal projection of $H$ to $H_n$. Define $u_n$ as the solution of the following stochastic differential equation; for each $v\in H_n$
\begin{equation}\label{eq:fs}
    (u_n(t) , v) = (u_n(0), v) + \int_0^t \langle F(u_n(s)) , v \rangle ds + \int_0^t \langle \sigma(u_n(s)) dW_n(s), v\rangle  
\end{equation}
with $u_n(0) = P_n u_0$.

\section{Existence and Uniqueness of Strong Solution}
In this section, we provide a priori estimates of\;\eqref{eq:fs}, establish the existence of a strong solution to\;\eqref{eq:af} by way of local monotonicity and later establish the uniqueness of the strong solution. The method of proving existence of strong solutions of SPDEs via monotonicity was introduced by Pardoux\;\cite{pard}. Proving the existence of a solution to the stochastic Navier-Stokes equations with additive noise by the method of local monotonicity was first shown by Menaldi and Sritharan\;\cite{men} and then by Sritharan and Sundar\;\cite{ss} for the case of the stochastic Navier-Stokes equations with multiplicative noise. This technique of proving existence of a strong solution by local monotonicity was later formalized by Liu and R\"{o}ckner\;\cite{lr}. For the sake of completion, we include the proof of existence and uniqueness of the strong solution of\;\eqref{eq:af}.

\begin{proposition}\label{prop:apriori}
    Let $\sigma(u)$ satisfy\;\Cref{ass:cv}. Let $u_n$ denote the unique strong solution of the finite system\;\eqref{eq:fs} in $C([0,T]; H_n)$, then, for $0\leq t\leq T$, the following estimates hold
    \begin{equation}\label{eq:eest}
        \mathbb{E}(|u_n(t)|^2) + \nu_0\int_0^t \mathbb{E}||u_n(s)||^2 ds + 2\nu_1\int_0^t\mathbb{E}||u_n(s)||^q_{V_{1,q}} ds \leq \widehat{C}\Big(1+\mathbb{E}|u_n(0)|^2  + \frac{1}{\nu_0}||f||_{L^2(0,T; V')}^2\Big)
    \end{equation}
    where $\widehat{C} = (1+ KTe^{KT})(1\vee KT)$ and
    \begin{equation}\label{eq:supeest}
    \begin{aligned}
        \mathbb{E}\sup_{0\leq t\leq T}\bigg(|u_n(t)|^2 + \nu_0\int_0^t ||u_n(s)||^2 ds + &2\nu_1\int_0^t||u_n(s)||^q_{V_{1,q}} ds\bigg)\\
        &\qquad\leq \widetilde{C}\Big(1+\mathbb{E}|u_n(0)|^2  + \frac{1}{\nu_0}||f||_{L^2(0,T; V')}^2\Big)
    \end{aligned}
    \end{equation}
    where $\widetilde{C} = (1 + 34KTe^{34KT} + 2e^{34KT})(1\vee 17KT)$.
\end{proposition}

\begin{proof}
    By the It\^{o} formula, we have, for $0\leq t \leq T$, 
    \begin{equation}\label{eq:ito}
    \begin{aligned}
        |u_n(t)|^2 + 2\nu_0 &\int_0^t \langle Au_n(s), u_n(s) \rangle ds + 2\nu_1\int_0^t \langle \mathcal{A}_1(u_n(s)), u_n(s) \rangle ds\\
        & = |u_n(0)|^2 + 2\int_0^t \langle f, u_n(s)\rangle ds + \int_0^t ||\sigma_n(u_n(s))||_{L_Q}^2 ds + 2\int_0^t \langle u_n(s), \sigma_n(u_n(s)) dW_n(s) \rangle.
    \end{aligned}
    \end{equation}
    Applying the properties of $A$ and $\mathcal{A}_1$ and the growth condition of $\sigma$, we get 
    \begin{align*}
        |u_n(t)|^2 + 2\nu_0 &\int_0^t ||u_n(s)||^2 ds + 2\nu_1 \int_0^t ||u_n(s)||_{V_{1,q}}^2 ds \\
        &\leq |u_n(0)|^2 + 2\int_0^t \langle f, u_n(s) \rangle ds + \int_0^t K(1 + |u_n(s)|^2) ds + M(t)
    \end{align*}
    where $M(t) := 2\int_0^t \langle u_n, \sigma_n(u_n) dW_n(s) \rangle$. By the Young inequality, we have 
    \begin{align*}
        |u_n(t)|^2 + \nu_0 &\int_0^t ||u_n(s)||^2 ds + 2\nu_1 \int_0^t ||u_n(s)||_{V_{1,q}}^2 ds \\
        &\leq |u_n(0)|^2 + \frac{1}{\nu_0} ||f||_{L^2(0,T; V')}^2 + Kt + K\int_0^t |u_n(s)|^2 ds + M(t).
    \end{align*}
    Now, define the stopping time 
    \[\tau_N := \inf\{t : |u_n(t)|^2 + \nu_0 \int_0^t||u_n(s)||^2 ds + 2\nu_1\int_0^t ||u_n(s)||_{V_{1,q}}^q ds > N\}\]
    where $\tau_N = T$ if the set is empty. Taking the expectation yields
    \begin{equation}\label{eq:stest}
    \begin{aligned}
        \mathbb{E}|u_{n}(t\wedge \tau_N)|^2 &+ \nu_0 \mathbb{E}\int_0^{t\wedge \tau_N} ||u_n(s)||^2 ds + 2\nu_1 \mathbb{E}\int_0^{t\wedge \tau_N} ||u_n(s)||_{V_{1,q}}^q ds\\
        &\leq \mathbb{E}|u_n(0)|^2 + \frac{1}{\nu_0} ||f||_{L^2(0,T; V')}^2 + KT + K\int_0^{t\wedge \tau_N} \mathbb{E}|u_n(s)|^2 ds.
    \end{aligned}
    \end{equation}
    Ignoring the second and third term on the left side of the inequality and applying the Gr\"{o}nwall inequality, we get 
    \begin{align*}
        \mathbb{E}|u_n(t\wedge \tau_N)|^2 \leq C_1\Big(1 + \mathbb{E}|u_n(0)|^2 + \frac{1}{\nu_0} ||f||_{L^2(0,T; V')}^2\Big)
    \end{align*}
    where $C_1 = e^{KT}(1\vee KT)$. Applying this estimate to the last term in\;\eqref{eq:stest}, we get 
    \begin{align*}
        \mathbb{E}|u_{n}(t\wedge \tau_N)|^2 + \nu_0 \mathbb{E}\int_0^{t\wedge \tau_N} ||u_n(s)||^2 ds + &2\nu_1 \mathbb{E}\int_0^{t\wedge \tau_N} ||u_n(s)||_{V_{1,q}}^q ds\\
        &\leq C_2\Big(1 + \mathbb{E}|u_n(0)|^2 + \frac{1}{\nu_0} ||f||_{L^2(0,T; V')}^2 \Big)
    \end{align*}
    where $C_2 = (1 + KT e^{KT})(1\vee KT)$. Taking $N\rightarrow \infty$, we arrive at the desired result\;\eqref{eq:eest}.\\

    Now, taking the supremum up to time $t\wedge \tau_N$ of\;\eqref{eq:ito} then taking the expectation gives us 
    \begin{equation}\label{eq:supeq}
    \begin{aligned}
        \mathbb{E}\sup_{0\leq s\leq t\wedge \tau_N} |u_n(s)|^2 &+ 2\nu_0\mathbb{E}\int_0^{t\wedge \tau_N} ||u_n(s)||^2 ds + 2\nu_1\mathbb{E}\int_0^{t\wedge \tau_N} ||u_n(s)||_{V_{1,q}}^q ds\\
        &= \mathbb{E}|u_n(0)|^2 + 2 \mathbb{E}\int_0^{t\wedge \tau_N} \langle f, u_n(s) \rangle ds + \mathbb{E}\int_0^{t\wedge \tau_N} ||\sigma_n(u_n(s))||_{L_Q}^2 ds\\
        &\quad + 2\mathbb{E}\sup_{0\leq s\leq t\wedge\tau_N} \int_0^{s} \langle u_n(r), \sigma_n(u_n(r)) dW_n(r)\rangle.
    \end{aligned}
    \end{equation}
    Applying the Burkholder-Davis-Gundy inequality to the martingale term, we have 
    \begin{equation}\label{eq:bdg}
    \begin{aligned}
        \mathbb{E}\sup_{0\leq s\leq t\wedge\tau_N} \bigg|\int_0^{s} \langle u_n(r), \sigma_n(u_n(r)) dW_n(r)\rangle\bigg| &\leq 2\sqrt{2} \mathbb{E}\bigg({\Big(\int_0^{t\wedge \tau_N}||\sigma_n(u_n)(s)||_{L_Q}^2 |u_n(s)|^2 ds\Big)}^{1/2}\bigg)\\
        &\leq 2\sqrt{2} \mathbb{E}\bigg(\sup_{0\leq s\leq t\wedge \tau_N} |u_n(s)|{\Big(\int_0^{t\wedge \tau_N} ||\sigma_n(u_n(s))||_{L_Q}^2\Big)}^{1/2}\bigg)\\
        &\leq \mathbb{E}\bigg(\frac{1}{4} \sup_{0\leq s\leq t\wedge \tau_N} |u_n(s)|^2 + 8\int_0^{t\wedge \tau_N} ||\sigma_n(u_n(s))||_{L_Q}^2 ds\bigg)\\
        &\leq \frac{1}{4}\mathbb{E}\sup_{0\leq s\leq t\wedge \tau_N} |u_n(s)|^2 + 8KT + 8K\mathbb{E}\int_0^{t\wedge \tau_N} |u_n(s)|^2 ds.
    \end{aligned}
    \end{equation}
    Substituting\;\eqref{eq:bdg} into\;\eqref{eq:supeq} gives 
    \begin{align*}
        \mathbb{E}\sup_{0\leq s\leq t\wedge \tau_N} |u_n(s)|^2 &+ 2\nu_0\mathbb{E}\int_0^{t\wedge \tau_N} ||u_n(s)||^2 ds + 2\nu_1\mathbb{E}\int_0^{t\wedge \tau_N} ||u_n(s)||_{V_{1,q}}^q ds\\
        &\leq \mathbb{E}|u_n(0)|^2 + \frac{1}{\nu_0}||f||_{L^2(0,T; V')}^2 +  \nu_0\mathbb{E}\int_0^{t\wedge \tau_N} ||u_n(s)||^2 ds\\
        &\quad + \frac{1}{2}\mathbb{E}\sup_{0\leq s\leq t\wedge \tau_N} |u_n(s)|^2 + 17KT + 17K\mathbb{E}\int_0^{t\wedge \tau_N} |u_n(s)|^2 ds.
    \end{align*}
    Through some simplification, we arrive at the following inequality
    \begin{equation}\label{eq:stsupest}
    \begin{aligned}
        \mathbb{E}\sup_{0\leq s\leq t\wedge \tau_N} |u_n(s)|^2 &+ 2\nu_0\mathbb{E}\int_0^{t\wedge \tau_N} ||u_n(s)||^2 ds + 4\nu_1\mathbb{E}\int_0^{t\wedge \tau_N} ||u_n(s)||_{V_{1,q}}^q ds\\
        &\leq 2\mathbb{E}|u_n(0)|^2 + \frac{2}{\nu_0}||f||_{L^2(0,T; V')}^2 + 34KT + 34K\int_0^{t\wedge \tau_N} \mathbb{E}\sup_{0\leq s\leq t\wedge \tau_N} |u_n(s)|^2 ds.
    \end{aligned}
    \end{equation}
    Ignoring the second and third term on the left side of the inequality and applying the Gr\"{o}nwall inequality, we get 
    \begin{equation}\label{eq:supee}
    \begin{aligned}
        \mathbb{E}\sup_{0\leq s\leq t\wedge \tau_N} |u_n(s)|^2 \leq C_3\Big(1 + \mathbb{E}|u_n(0)|^2 + \frac{1}{\nu_0}||f||_{L^2(0,T; V')}^2\Big)
    \end{aligned}
    \end{equation}
    where $C_3 = 2e^{34KT}(1\vee 17KT)$. Substituting this inequality into the last term of\;\eqref{eq:stsupest}, we get 
    \begin{equation}\label{eq:stsupest2}
    \begin{aligned}
        \nu_0\mathbb{E}&\int_0^{t\wedge \tau_N} ||u_n(s)||^2 ds + 2\nu_1\mathbb{E}\int_0^{t\wedge \tau_N} ||u_n(s)||_{V_{1,q}}^q ds\\
        &\leq \mathbb{E}|u_n(0)|^2 + \frac{1}{\nu_0}\int_0^{T} ||f||_{V'}^2 ds + 17KT + 17KTC_3\Big(1+ \mathbb{E}|u_n(0)|^2 + \frac{1}{\nu_0}||f||_{L^2(0,T; V')}^2\Big)\\
        &\leq C_4\Big(1 + \mathbb{E}|u_n(0)|^2 + \frac{1}{\nu_0}||f||_{L^2(0,T; V')}^2\Big)
    \end{aligned}
    \end{equation}
    where $C_4 = (1+34KTe^{34KT})(1\vee 17KT)$. Combining\;\eqref{eq:supee} and\;\eqref{eq:stsupest2} gives 
    \begin{align*}
        \mathbb{E}\sup_{0\leq s\leq t\wedge \tau_N} |u_n(s)|^2 + \nu_0\mathbb{E}\int_0^{t\wedge \tau_N} ||u_n(s)||^2 ds &+ 2\nu_1\mathbb{E}\int_0^{t\wedge \tau_N} ||u_n(s)||_{V_{1,q}}^q ds\\
        &\leq C_5\Big(1 + \mathbb{E}|u_n(0)|^2 + \frac{1}{\nu_0}||f||_{L^2(0,T; V')}^2\Big)
    \end{align*}
    where $C_5 = (1 + 34KTe^{KT} + 2e^{34KT})(1\vee 17KT)$. Taking $N\rightarrow \infty$ gives the desired result.
\end{proof}

Now, that we have the a priori estimates, we are ready to state the theorem on existence and uniqueness of the strong solution.

\begin{theorem}\label{thm:solution}
    Let $\sigma$ satisfy\;\Cref{ass:cv}. If $2L<\nu_0$ then, for each initial condition $u_0\in H$, there exists a unique strong solution $u$ of\;\eqref{eq:af} in the Lusin space
    \[L^2(\Omega; C(0,T; H))\cap L^q(\Omega; L^q(0,T; V_{1,q})).\]
\end{theorem}

\begin{proof}
    Let $F(u) := -\nu_0 Au - \nu_1\mathcal{A}_1(u) - B(u)$. Then we have 
    \[u_n(t) = u_n(0) + \int_0^t F(u_n(s)) ds + \int_0^t f\ ds + \int_0^t \sigma_n(u_n(s))dW(s).\]
    By the a priori estimates and the Banach-Alaoglu theorem, we have the following limits along a subsequence of the Galerkin approximations
    \begin{align*}
        &u_n \rightarrow u \text{ weakly in } L^q(\Omega; L^q(0,T; V_{1,q}))\\
        &u_n \rightarrow u \text{ weak star in } L^2(\Omega; L^\infty(0,T; H)).
    \end{align*}
    Since $F(u_n)$ is bounded in $L^{q'}(\Omega; L^{q'}(0,T; V_{1,q}'))$
    \[F(u_n) \rightarrow F_0 \text{ weakly in } L^{q'}(\Omega; L^{q'}(0,T; V_{1,q}')).\]
    Here $q'$ denotes the H\"{o}lder conjugate of $q$ meaning $q' = \frac{q}{q-1}$. By the assumptions on $\sigma$,
    \[\sigma_n(u_n) \rightarrow S \text{ weakly in } L^2(\Omega;L^2(0,T; L_Q)).\]
    It follows that $u$ satisfies 
    \[u(t) = u_0 + \int_0^t F_0\ ds + \int_0^t f\ ds + \int_0^t S\ dW(s).\]
    For any adapted process $v$ in $L^\infty (\Omega \times [0,T], H_m)$ with $m\leq n$, we define 
    \[r(t,\omega) := \frac{1}{\nu_0^3} \int_0^t ||v(s)||_{L^4}^4 ds\]
    as an adapted, continuous real valued process. By monotonicity,
    \begin{align*}
        X_n(t):=& \mathbb{E}\int_0^t e^{-r(s)} \langle F(u_n(s)) - F(v(s)), u_n(s) - v(s) \rangle ds\\
        & - \mathbb{E}\int_0^t \frac{dr(s)}{ds} e^{-r(s)} |u_n(s) - v(s)|^2 ds\\
        & + \mathbb{E} \int_0^t e^{-r(s)} ||\sigma_n(u_n(s)) - \sigma_n(v(s))||^2_{L_Q}ds \leq 0.
    \end{align*}
    We will decompose $X_n$ into the sum $Y_n + Z_n$ where
    \begin{align*}
        Y_n := &\mathbb{E}\int_0^t e^{-r(s)} \langle F(u_n(s)), u_n(s) \rangle ds\\
        & - \mathbb{E} \int_0^t \frac{dr(s)}{ds}e^{-r(s)}|u_n(s)|^2 ds + \mathbb{E}\int_0^t e^{-r(s)} ||\sigma_n(u_n(s))||^2_{L_Q} ds
    \end{align*}
    and
    \begin{align*}
        Z_n :=& \mathbb{E}\int_0^t e^{-r(s)}\Big(\langle F(v(s)), v(s) \rangle - \langle F(v(s)) , u_n(s) \rangle - \langle F(u_n(s)) , v(s) \rangle\Big)ds\\
        & - \mathbb{E}\int_0^t \frac{dr(s)}{ds}e^{-r(s)}\Big(|v(s)|^2 - 2(u_n(s),v(s))\Big) ds\\
        & + \mathbb{E}\int_0^t e^{-r(s)}\Big(||\sigma_n(v(s))||^2_{L_Q} - 2\big(\sigma_n(u_n(s)), \sigma_n(v(s))\big)\Big) ds.
    \end{align*}
    By the It\^o formula, we have
    \[Y_n = \mathbb{E}\big(e^{-r(t)} |u_n(t)|^2 - |u_n(0)|^2\big).\]
    It follows that
    \begin{align*}
        \liminf_{n\rightarrow\infty} Y_n &\geq \mathbb{E}\big(e^{-r(t)} |u(t)|^2 - |u_0|^2\big)\\
        & = \mathbb{E}\int_0^t e^{-r(s)} \langle F_0, u(s)\rangle ds - \mathbb{E}\int_0^t \frac{dr(s)}{ds} e^{-r(s)} |u(s)|^2 ds\\
        &\phantom{= \;} + \mathbb{E}\int_0^t e^{-r(s)} ||S(s)||^2_{L_Q} ds + \mathbb{E}\int_0^t e^{-r(s)}\langle S(s) dW(s), u(s)\rangle.
    \end{align*}
    Similarly, we have 
    \begin{align*}
        \liminf_{n\rightarrow \infty} Z_n &\geq \mathbb{E}\int_0^t e^{-r(s)}\Big(\langle F(v(s)), v(s) \rangle - \langle F(v(s)) , u(s) \rangle - \langle F_0 , v(s) \rangle\Big)ds\\
        & - \mathbb{E}\int_0^t \frac{dr(s)}{ds}e^{-r(s)}\Big(|v(s)|^2 - 2(u(s),v(s))\Big) ds\\
        & + \mathbb{E}\int_0^t e^{-r(s)} \Big(||\sigma(v(s))||^2_{L_Q} - 2\big(S(s), \sigma(v(s))\big)\Big) ds.
    \end{align*}
    Combining these two inequalities, we arrive at 
    \begin{equation}\label{eq:ineq}
        \begin{aligned}
            \mathbb{E}\int_0^t e^{-r(s)}\langle F_0 - F(v(s)), u(s)-v(s) \rangle ds &- \mathbb{E}\int_0^t \frac{dr(s)}{ds}e^{-r(s)} |u(s)-v(s)|^2 ds\\
        &\qquad + \mathbb{E}\int_0^t e^{-r(s)}||S(s) - \sigma(v(s))||^2_{L_Q} ds \leq \liminf_{n\rightarrow \infty} X_n \leq 0.    
        \end{aligned}
    \end{equation}
    {Note that this inequality holds for every $v\in L^\infty(\Omega \times [0,T], H_m)$ and any $m$. It is immediate that the inequality holds for any
    \[v\in \bigcup_{m=1}^\infty L^\infty(\Omega \times [0,T], H_m) = L^\infty(\Omega \times [0,T], \bigcup_{m=1}^\infty H_m).\]
    Since $\bigcup_{m} H_m$ is dense in $H$, a standard approximation argument shows that
    \[L^\infty(\Omega \times [0,T], \bigcup_{m=1}^\infty H_m) \text{ is dense in } L^2(\Omega, L^\infty(0,T; H)).\]
    It follows that inequality\;\eqref{eq:ineq} holds for any adapted process $v$ in $L^2(\Omega, L^\infty(0,T; H))$. So, we take $v:= u$.}
    Then 
    \begin{align*}
        \mathbb{E}\int_0^T e^{-r(t)} ||S(t) - \sigma(u(s))||^2_{L_Q} ds = 0.
    \end{align*}
    Hence $S(s) = \sigma(u(s))$. Now, take $v:= u - \lambda w$ where $w$ is an arbitrary adapted process in $L^2(\Omega, L^\infty(0,T; H))$. Then
    \begin{align*}
        &\mathbb{E}\int_0^t e^{-r(s)}\langle F_0 - F(u(s)-\lambda w(s)), \lambda w(s) \rangle ds - \mathbb{E}\int_0^t \frac{dr(s)}{ds}e^{-r(s)} |\lambda w(s)|^2 ds\\
        &\phantom{\qquad} = \lambda \mathbb{E}\int_0^t e^{-r(s)}\langle F_0 - F(u(s)-\lambda w(s)) , w(s)\rangle ds - \lambda^2 \mathbb{E}\int_0^t \frac{dr(s)}{ds} e^{-r(s)} |w(s)|^2 ds\\
        &\phantom{\qquad} \leq 0.
    \end{align*}
    Dividing by $\lambda$ then taking $\lambda \rightarrow 0$, we obtain
    \[\mathbb{E}\int_0^t e^{-r(s)} \langle F_0 - F(u(s)) , w(s) \rangle ds \leq 0.\]
    Since $w$ is arbitrary, $F_0$ can be identified with $F(u)$. Thus a strong solution exists.\\

    Suppose $u,v$ are two solutions to\;\eqref{eq:af} with initial conditions $u_0$ and $v_0$ respectively. For $N>0$, define $\tau^u_N:= \inf\{t\leq T : |u(t)|\geq N\}$ and $\tau^v_N := \inf\{t\leq T : |v(t)|\geq N\}$. Let $\tau_N = \min\{\tau^u_N, \tau^v_N\}$. By the Markov inequality,
    \begin{align*}
        \mathbb{P}\{\tau_N < T\}\leq \mathbb{P}\{\max\{\sup_{t\in[0,T]}|u(t)|,\sup_{t\in[0,T]}|v(t)|\}\geq N\}\leq \frac{C}{N^2}.
    \end{align*}
    Thus, as $N\rightarrow \infty$, $\tau_N \rightarrow T$ a.s. Define $w(t)= u(t)-v(t)$ and $w_0 = u_0 - v_0$. Then
    \begin{align*}
    w(t) = w_0 -\nu_0 \int_0^t Aw(s)\ ds\ -\ &\nu_1\int_0^t (\mathcal{A}_1(u(s)) - \mathcal{A}_1(v(s)))\ ds\\
    & - \int_0^t (B(u(s)) - B(v(s)))\ ds + \int_0^t (\sigma(u(s)) - \sigma(v(s)))\ dW(s).
    \end{align*}
    Applying the It\^{o} formula to $w\mapsto |w(t\wedge \tau_N)|^2$, we arrive at 
     \begin{align*}
        d|w(t)|^2 + (2\nu_0||w(t)||^2 + &2\nu_1\langle \mathcal{A}_1(u(t)) - \mathcal{A}_1(v(t)), w(t)\rangle + 2\langle B(u(t)) - B(v(t)), w(t)\rangle)dt\\
        &=||\sigma(u(t))-\sigma(v(t))||^2_{L_Q}dt + 2\langle w(t), (\sigma(u(t)) - \sigma(v(t)))dW(t) \rangle
    \end{align*}
    By\;\eqref{eq:nl},\;\eqref{eq:bl}, and the assumptions on $\sigma$, we have
    \begin{align*}
        d|w(t)|^2 + (2\nu_0 ||w(t)||^2 + 2\nu_1\gamma||w(t)||^q_{V_{1,q}} &+ 2\sqrt{2}|w(t)|\ ||w(t)||\ ||u(t)||)dt\\
        &\leq L|w(t)|^2 dt + 2\langle w(t), (\sigma(u(t)) - \sigma(v(t)))dW(t) \rangle.
    \end{align*}
    Using the Poincar\'{e} and Young inequalities, we get
    \begin{align*}
        d|w(t)|^2 + (\nu_0\lambda_1 |w(t)|^2 &+ 2\nu_1\gamma||w(t)||^q_{V_{1,q}})dt\\
        &\leq \bigg(\frac{2}{\nu_0}||u(t)||^2+ L\bigg)|w(t)|^2 dt + 2\langle w(t), (\sigma(u(t)) - \sigma(v(t)))dW(t) \rangle
    \end{align*}
    Rearranging some terms and removing the positive term $2\nu_1\gamma||w(t)||_{V_{1,q}}^q$ allows us to write 
    \[d|w(t)|^2 + \bigg(\nu_0\lambda_1 -\frac{2}{\nu_0}||u(t)||^2- L\bigg)|w(t)|^2 dt \leq 2\langle w(t), (\sigma(u(t)) - \sigma(v(t)))dW(t) \rangle.\]
    Let
    \begin{equation}\label{eq:gamma}
       \Gamma(t) := (\nu_0\lambda_1 - L)t - \frac{2}{\nu_0}\int_0^t ||u(s)||^2 ds.  
    \end{equation}
    Then multiplying by $e^{\Gamma(t)}$ and integrating up to $t\wedge \tau_N$ gives 
    \[e^{\Gamma(t\wedge \tau_N)}|w(t\wedge \tau_N)|^2 - |w_0|^2 \leq 2\int_0^{t\wedge \tau_N} e^{\Gamma(s)}\langle w(s), (\sigma(u(s)) - \sigma(v(s)))dW(s) \rangle.\] 
    Observe that, for a fixed $N$, the process
    \[M_N(t) := 2\int_0^{t\wedge \tau_N} e^{\Gamma(s)} \langle w(s), (\sigma(u(s)) - \sigma(v(s)))dW(s) \rangle\]
    is a martingale since 
    \[\mathbb{E}\int_0^{t\wedge \tau_N} |w(s)|^2 ||\sigma(u(s)) - \sigma(v(s))||_{L_Q}^2 ds \leq L\mathbb{E}\int_0^{t\wedge \tau_N} |w(s)|^4 ds \leq 16LTN^4 < \infty.\]
    Taking the expectation gives 
    \begin{equation}\label{eq:gronwall}
    \begin{aligned}
        \mathbb{E}e^{\Gamma(t\wedge \tau_N)}|w(t\wedge \tau_N)|^2 \leq \mathbb{E}|w_0|^2. 
    \end{aligned}
    \end{equation}
    Suppose $u$ and $v$ have the same initial condition, i.e. $u_0 = v_0$. Then $w_0 = 0$. Hence\;\eqref{eq:gronwall} implies $w(t\wedge \tau_N) = 0$ almost surely. Thus, taking $N\rightarrow \infty$, we conclude that $u=v$ for all $t\in [0,T]$ almost surely.
\end{proof}

\section{Existence of Ergodic Invariant Measure}
For this section, we assume the initial condition $u_0 \in H$ is deterministic and the forcing $f$ is time independent. Define the transition probability measure
\[P_t(x,B) = P\{u_t^x \in B\}\]
for all $t\geq 0$, $x\in H$, and $B\in \mathcal{B}(H)$ where $u_t^x$ denotes the solution of\;\eqref{eq:af} starting at $x$. Then $P_t (\cdot,B)$ is an $\mathcal{B}(H)$-measurable function for all fixed $t$ and $B$. The Markov property yields the Chapman-Kolmogorov equation
\[P_{t+s}(x,B) = \int_H P_t(x,dy) P_s(y,B).\]
We define the semigroup of linear operators $P_t$ on $B_b(H)$ as follows 
\[P_t\phi(x) =\mathbb{E}\phi(u_t^x) = \int_H \phi(y) P_t(x,dy).\]
The semigroup $P_t$ is strongly continuous and is called a Markov semigroup since $P_t$ is a positivity preserving contraction and $P_t1=1$. Now, we define invariant measure corresponding to the transition probability function $P_t$ and establish a definition for an ergodic invariant measure.

\begin{defn}
    A probability measure $\mu$ on $(H,\mathcal{B}(H))$ is called an invariant measure for the transition probability function $P_t$ if it satisfies 
    \[\mu(B) = \int_H P_t(x,B) d\mu(x)\]
    for all $B\in \mathcal{B}(H)$, or equivalently, if for all $\phi \in C_b(H)$ and all $t\geq 0$,
    \[\int_H \phi(x) d\mu(x) = \int_H (P_t \phi)(x) d\mu(x).\]
\end{defn}

\begin{theorem}\label{thm:inv}
    Let\;\Cref{ass:cv} and the inequality 
    \begin{equation}
        \frac{\nu_0\lambda_1}{2} + 2\nu_1\tilde{C} > K
    \end{equation}
    hold. Then, the transition probability function $P_t$ corresponding to\;\eqref{eq:af} admits at least one invariant measure.
\end{theorem}

\begin{proof}
    Existence of an invariant measure will follow from Krylov-Bogolyubov Theorem provided we show that the probability transition function $P_t$ is Feller and the tightness property holds. To show $P_t$ is Feller, let $u$ and $u_n$ denote the solutions to\;\eqref{eq:af} with initial conditions $u_0, u_{0,n}$, respectively, where $u_{0,n}$ converges to $u_0$ a.s. Now, observe that, by\;\eqref{eq:gronwall}, 
    \[\mathbb{E}[e^{\Gamma(t)} |u(t) - u_n(t)|^2] \leq |u_0 - u_{0,n}|^2.\]
    Taking $n\rightarrow \infty$, it follows that $u_n \rightarrow u$ a.s. By the continuity of $\phi$, we have that $\phi(u_n(t)) \rightarrow \phi(u(t))$ a.s. Then, by dominated convergence theorem, $\mathbb{E}\big(\phi(u_n(t))\big)\rightarrow \mathbb{E}\big(\phi(u(t))\big)$. Thus, as $n\rightarrow \infty$, 
    \[|P_t\phi(u_{0,n}) - P_t \phi(u_0)|\leq \mathbb{E}|\phi(u_n(t)) - \phi(u(t))| \rightarrow 0\]
    and hence $P_t$ is Feller. Now, we show the tightness property holds. For the purposes of simplicity, let $u_0 = 0$. We need to show that the family of probability measures $\mu_t$, defined by
    \[\mu_t(B) := \frac{1}{t} \int_0^t P_s(0,B) ds\]
    for all $B \in \mathcal{B}(H)$ and $t>0$ is tight in $H$. Let $B_n$ denote the closed ball in $V_{1,q}$ of radius $n$. Note that $B_n$ is a compact subset of $H$ since the embedding $V_{1,q} \hookrightarrow H$ is continuous and compact. By the Markov inequality and\;\eqref{eq:eest}, we have 
    \begin{align*}
        \lim_{n\rightarrow \infty}\sup_{t} \mu_t(B_n^c) &= \lim_{n\rightarrow \infty}\sup_{t} \frac{1}{t}\int_0^t \mathbb{P}\{||u(s)||_{V_{1,q}} \geq n\} ds\\
        &\leq \lim_{n\rightarrow \infty}\sup_{t} \frac{1}{tn^q} \int_0^t \mathbb{E}||u(s)||_{V_{1,q}}^q ds\\
        &\leq \lim_{n\rightarrow \infty} \frac{C}{n^q}\\
        &=0
    \end{align*}
    where $C$ is a suitable constant. It follows that $\mu_t$ is tight in $H$. Thus, by Krylov-Bogolyubov Theorem, the transition probability function $P_t$ admits at least one invariant measure.
\end{proof}

\begin{defn}
    Let $\mu$ be an invariant measure for $P_t$. A set $A \in \mathcal{B}(H)$ is called an invariant set with respect to the semigroup $P_t$ if for all $t\geq 0$
    \[P_t \mathbbm{1}_A = \mathbbm{1}_A \quad \mu\text{-a.s.}\]
    The measure $\mu$ is called ergodic if $\mu(A)$ is either 0 or 1 for any invariant set $A$.
\end{defn}

Let $\mathcal{I}$ denote the set of invariant probability measures for the semigroup $P_t, t\geq 0$. Since Theorem 4.1 proves the existence of at least one invariant measure, $\mathcal{I}$ is nonempty. The following proposition is due to\;\cite{dpz}.

\begin{proposition}[\cite{dpz}, Proposition 3.2.7]\label{prop:dpz}
    An invariant probability measure for the semigroup $P_t, t\geq 0$, is ergodic if and only if it is an extremal point of $\mathcal{I}$.
\end{proposition}

{Since\;\Cref{thm:inv} shows $\mathcal{I}$ is nonempty, we are left to show that $\mathcal{I}$ contains at least one extremal point. This will be accomplished through the Krein-Milman Theorem.}
{
\begin{theorem}[Krein-Milman]
    Every nonempty compact convex subset of a Hausdorff locally convex topological vector space has at least one extreme point.
\end{theorem}}

Before stating the existence of an extremal point, we need to determine estimates on the moments of the solution.

\begin{proposition}[Moment Estimates]\label{prop:moment}
    Let $\sigma$ satisfy\;\Cref{ass:cv} and assume
    \[\frac{\nu_0\lambda_1}{2} + 2\nu_1\tilde{C} > K.\]
    For $p\in \left[2, 1+ \frac{\nu_0\lambda_1 + 4\tilde{C}\nu_1}{2K}\right)$, the solution $u$ of the NSNE with shear-thickening viscosity satisfies 
    \begin{equation}\label{eq:moment}
        \mathbb{E}[|u(t)|^p] \leq |u_0|^p e^{-at} + \frac{C}{a} \qquad \forall t> 0  
    \end{equation}
    where 
    \[a = \begin{cases}
     \frac{\nu_0\lambda_1}{2} + 2\nu_1 \tilde{C} - K & \text{if } p = 2\\
    \frac{p\nu_0\lambda_1}{4} + \tilde{C}p\nu_1 - \frac{Kp(p-1)}{2} \quad & \text{if } p \in (2, 1+ \frac{\nu_0\lambda_1 + 4\tilde{C}\nu_1}{2K})
    \end{cases}\]
    and 
    \[C = \begin{cases}
        \frac{2}{3\nu_0}||f||_{V'}^2 + K + 2\nu_1\tilde{C} & \text{if } p = 2\\
        \frac{2^{\frac{p}{2}-1}}{{\nu_0}^{p-1}{\lambda_1}^{\frac{p}{2}-1}}||f||^p_{V'} + \frac{{(p-1)}^{\frac{p}{2}} {(p-2)}^{\frac{p}{2}-1}}{{(\nu_0\lambda_1)}^{\frac{p}{2}-1}}K^{\frac{p}{2}} + \tilde{C}p\nu_1 \quad & \text{if } p \in (2, 1+ \frac{\nu_0\lambda_1 + 4\tilde{C}\nu_1}{2K})
    \end{cases}\]
    Here $\tilde{C} = \lambda_1^{q/2} |G|^{-\frac{q}{2} + 1}$ and $|G|$ is the size of the bounded domain.
\end{proposition}

\begin{proof}
    For \(p = 2\), we obtain the following equality from the It\^{o} formula
    \begin{equation}\label{eq:eq1}
    \begin{aligned}
        |u(t)|^2 + 2\nu_0\int_0^t ||u(s)||^2 ds &+ 2 \nu_1\int_0^t ||u(s)||^q_{V_{1,q}} ds\\
        & = |u_0|^2 + 2\int_0^t \langle f,u(s)\rangle ds + \int_0^t ||\sigma(u(s))||^2_{L_Q} ds + 2\int_0^t\langle u(s), \sigma(u(s)) dW(s)\rangle.
    \end{aligned}
    \end{equation}
    By the Young inequality and Assumption\;\ref{M2}, we get 
    \begin{equation}\label{eq:est1}
    \begin{aligned}
        |u(t)|^2 + \frac{\nu_0}{2}\int_0^t ||u(s)||^2 ds &+ 2 \nu_1\int_0^t ||u(s)||^q_{V_{1,q}} ds\\
        & \leq |u_0|^2 + \frac{2}{3\nu_0}\int_0^t ||f||_{V'}^2 + Kt + K\int_0^t |u(s)|^2 ds + 2\int_0^t\langle u(s), \sigma(u(s)) dW(s)\rangle.
    \end{aligned}
    \end{equation}
    Taking the expectation of\;\eqref{eq:est1} gives 
    \begin{equation}
    \begin{aligned}
        \mathbb{E}|u(t)|^2 + \frac{\nu_0}{2} \int_0^t \mathbb{E}||u(s)||^2 ds &+ 2\nu_1\int_0^t \mathbb{E}||u(s)||_{V_{1,q}}^q ds \\
        &\leq |u_0|^2 + \frac{2}{3\nu_0}||f||_{V'}^2t + Kt + K\int_0^t \mathbb{E}|u(s)|^2 ds.
    \end{aligned}
    \end{equation}
    By using the Poincar\'{e} inequality and\;\Cref{prop:lowerbound} with $p= 2$, we get 
    \begin{align*}
        \mathbb{E}|u(t)|^2 + \Big(\frac{\nu_0\lambda_1}{2} + 2\nu_1\tilde{C} - K\Big)\int_0^t \mathbb{E}|u(s)|^2 ds \leq |u_0|^2 + \Big(\frac{2}{3\nu_0}||f||_{V'}^2 + K + 2\nu_1\tilde{C}\Big)t.
    \end{align*} 
    Applying the Gr\"{o}nwall inequality, we arrive at 
    \[\mathbb{E}|u(t)|^2 \leq |u_0|^2 e^{-at} + \frac{C}{a}\]
    where $a =  \frac{\nu_0\lambda_1}{2} + 2\nu_1 \tilde{C} - K$ and $C = \frac{2}{3\nu_0}||f||_{V'}^2 + K + 2\nu_1\tilde{C}$.\\

    Let $p\in \left(2, 1+ \frac{\nu_0\lambda_1 + 4\tilde{C}\nu_1}{2K}\right)$. Then, by using the It\^{o} formula in differential form, 
    \begin{equation}\label{eq:eqp}
    \begin{aligned}
        d|u(t)|^p +& (p\nu_0 |u(t)|^{p-2}||u(t)||^2 + p\nu_1|u(t)|^{p-2}||u(t)||_{V_{1,q}}^q)dt\\
        & \leq p|u(t)|^{p-2}\langle f,u(t) \rangle dt + \frac{p(p-1)}{2} |u(t)|^{p-2} ||\sigma(u(t))||_{L_Q}^2 dt + p|u(t)|^{p-2}\langle u, \sigma dW(t)\rangle.
    \end{aligned}
    \end{equation}
    Now, we determine upper bounds for the first and second term on the right side of the equality. The following upper bounds are achieved by applying the basic inequality, due to Young, multiple times
    \begin{align*}
        p|u|^{p-2}\langle f,u \rangle &\leq p|u|^{p-2}||u||\, ||f||_{V'}\\
        &\leq p|u|^{p-2}\bigg( \frac{\nu_0}{2} ||u||^2 + \frac{1}{2\nu_0}||f||^2_{V'}\bigg)\\
        &= \frac{p\nu_0}{2}|u|^{p-2}||u||^2 + \frac{p}{2\nu_0} |u|^{p-2}||f||^2_{V'}\\
        &\leq \frac{p\nu_0}{2}|u|^{p-2}||u||^2 + \frac{1}{2\nu_0}\bigg(\frac{\nu_0^2\lambda_1(p-2)}{2}|u|^p + \frac{2^{\frac{p}{2}}}{{(\nu_0^2 \lambda_1)}^{\frac{p}{2} - 1}}||f||^p_{V'}\bigg)\\
        &= \frac{p\nu_0}{2}|u|^{p-2}||u||^2 + \frac{\nu_0\lambda_1(p-2)}{4}|u|^p + \frac{2^{\frac{p}{2}-1}}{{\nu_0}^{p-1}{\lambda_1}^{\frac{p}{2}-1}}||f||^p_{V'}
    \end{align*} 
    \begin{align*}
        \frac{p(p-1)}{2} |u|^{p-2} ||\sigma(u)||_{L_Q}^2 &\leq \frac{p(p-1)}{2}|u|^{p-2}K(1+|u|^2)\\
        &= \frac{p(p-1)}{2}|u|^{p-2}K + \frac{Kp(p-1)}{2}|u|^{p}\\
        &\leq \frac{p(p-1)}{2}\bigg(\frac{\nu_0\lambda_1}{p(p-1)}|u|^p + \frac{2}{p}{\Big(\frac{(p-1)(p-2)}{\nu_0\lambda_1}\Big)}^{\frac{p}{2}-1}K^{\frac{p}{2}}\bigg)\\
        &\qquad + \frac{Kp(p-1)}{2}|u|^{p}\\
        &= \bigg(\frac{\nu_0\lambda_1}{2} + \frac{Kp(p-1)}{2}\bigg)|u|^p + \frac{{(p-1)}^{\frac{p}{2}} {(p-2)}^{\frac{p}{2}-1}}{{(\nu_0\lambda_1)}^{\frac{p}{2}-1}}K^{\frac{p}{2}}.
    \end{align*}
    Applying these inequalities and\;\Cref{prop:lowerbound}, we get 
    \begin{align*}
        d|u(t)|^p &+ (p\nu_0 |u(t)|^{p-2}||u(t)||^2 + \tilde{C}p\nu_1(|u(t)|^p -1))dt\\
        &\leq \bigg(\frac{p\nu_0}{2}|u(t)|^{p-2}||u(t)||^2 + \bigg(\frac{\nu_0\lambda_1(p-2)}{4} + \frac{\nu_0\lambda_1}{2} + \frac{Kp(p-1)}{2}\bigg)|u(t)|^p + \frac{2^{\frac{p}{2}-1}}{{\nu_0}^{p-1}{\lambda_1}^{\frac{p}{2}-1}}||f||^p_{V'}\\
        &\qquad + \frac{{(p-1)}^{\frac{p}{2}} {(p-2)}^{\frac{p}{2}-1}}{{(\nu_0\lambda_1)}^{\frac{p}{2}-1}}K^{\frac{p}{2}}\bigg)dt + p|u(t)|^{p-2}\langle u(t), \sigma(u(t)) dW(t)\rangle.
    \end{align*}
    Simplifying gives us 
    \begin{align*}
        d|u(t)|^p &+ \bigg(\frac{p\nu_0\lambda_1}{4} + \tilde{C}p\nu_1 - \frac{Kp(p-1)}{2}\bigg)|u(t)|^p dt\\
        &\leq \bigg(\frac{2^{\frac{p}{2}-1}}{{\nu_0}^{p-1}{\lambda_1}^{\frac{p}{2}-1}}||f||^p_{V'} + \frac{{(p-1)}^{\frac{p}{2}} {(p-2)}^{\frac{p}{2}-1}}{{(\nu_0\lambda_1)}^{\frac{p}{2}-1}}K^{\frac{p}{2}} + \tilde{C}p\nu_1\bigg)dt + p|u(t)|^{p-2}\langle u(t), \sigma(u(t)) dW(t)\rangle.
    \end{align*}
    Taking the expectation, we get 
    \[\frac{d}{dt} \mathbb{E}[|u(t)|^p] \leq -a\mathbb{E}[|u(t)|^p] + C\]
    where 
    \begin{align*}
        a &= \frac{p\nu_0\lambda_1}{4} + \tilde{C}p\nu_1 - \frac{Kp(p-1)}{2}\\
        C &= \frac{2^{\frac{p}{2}-1}}{{\nu_0}^{p-1}{\lambda_1}^{\frac{p}{2}-1}}||f||^p_{V'} + \frac{{(p-1)}^{\frac{p}{2}} {(p-2)}^{\frac{p}{2}-1}}{{(\nu_0\lambda_1)}^{\frac{p}{2}-1}}K^{\frac{p}{2}} + \tilde{C}p\nu_1.
    \end{align*}
    By the Gr\"{o}nwall inequality, we obtain 
    \[\mathbb{E}[|u(t)|^p] \leq |u_0|^p e^{-at} + \frac{C}{a}\]
    as desired.
\end{proof}

Now, we are ready for the existence of an extremal point of $\mathcal{I}$.

\begin{proposition}\label{prop:extremal}
    Let $\sigma$ satisfy\;\Cref{ass:cv} and assume
    \[\frac{\nu_0\lambda_1}{2} + 2\nu_1\tilde{C} > K.\]
    The set of invariant measures, $\mathcal{I}$, for the probability transition function $P_t$ has an extremal point.
\end{proposition}
{
\begin{proof}
    Consider the set of finitely additive signed measure on $H$, ${(C_b(H))}^*$. The set ${(C_b(H))}^*$ is a Hausdorff locally convex topological vector space under the weak-star topology such that $\mathcal{I} \subset {(C_b(H))}^*$. By Theorem 5.1, we have that $\mathcal{I}$ is nonempty. It remains to be shown that $\mathcal{I}$ is closed, compact, and convex. We begin with convexity. To show that $\mathcal{I}$ is convex, we must show that the measure
    \[\mu = (1-\lambda)\mu_1 + \lambda\mu_2\]
    is invariant where $\mu_1,\mu_2 \in \mathcal{I}$ and $\lambda\in [0,1]$. Observe that 
    \begin{align*}
        \int_H (P_t\phi) d\mu &= (1-\lambda)\int_H (P_t\phi) d\mu_1 + \lambda\int_H (P_t \phi) d\mu_2\\
        &= (1-\lambda)\int_H \phi\ d\mu_1 + \lambda\int_H \phi\  d\mu_2\\
        &= \int_H \phi \ d\mu.
    \end{align*}
    Thus, $\mathcal{I}$ is convex. Now, we move on to closedness. Suppose $\mu_n \Rightarrow \mu$ where $\mu_n \in \mathcal{I}$ for each $n$. To show that $\mathcal{I}$ is closed, we need to show that $\mu \in \mathcal{I}$. By the invariance and convergence of $\mu_n$, we can say
    \[\int_H (P_t \phi) d\mu_n = \int_H \phi \ d\mu_n \rightarrow \int_H \phi \ d\mu.\]
    By the Feller property of $P_t$, we have that $P_t \phi \in C_b(H)$, so we can pass the limit such that 
    \[\int_H (P_t \phi) d\mu_n \rightarrow \int_H(P_t \phi) d\mu.\]
    By uniqueness of limits, 
    \[\int_H \phi\ d\mu = \int_H(P_t \phi) d\mu.\]
    Hence, $\mathcal{I}$ is closed. To finish, we now show compactness of $\mathcal{I}$. Consider the set
    \[B_n = \{x \in H : |x| \leq n\}.\]
    Note that $B_n$ is a compact subset of $H$. For an invariant measure, $\mu$, we have by the Markov inequality 
    \begin{align*}
        \mu(B_n^c) &= \mu\{|x| > n\}\\ 
        &\leq \frac{1}{n^p} \int_H |x|^p d\mu(x)
    \end{align*}
    for $p\in \left[2, 1+ \frac{\nu_0\lambda_1 + 4\tilde{C}\nu_1}{2K}\right)$. By the invariance of $\mu$ and the moment estimate\;\eqref{eq:moment}, we have
    \begin{align*}
        \mu(B_n^c) &\leq \frac{1}{n^p} \int_H \mathbb{E}_x |u(t)|^p d\mu(x)\\
        &\leq \frac{1}{n^p}\int_H \Big(|x|^p e^{-at} + \frac{C}{a}\Big) d\mu(x)\\
        &= \frac{1}{n^p} e^{-at} \int_H |x|^p d\mu(x) + \frac{C}{n^p a}
    \end{align*}
    Taking $t \rightarrow \infty$ gives $\mu(B_n^c) \leq \frac{C}{n^p a}$. Let $\varepsilon >0$. Choose $n > \sqrt[p]{\frac{C}{a\varepsilon}}$. Then $\mu(B_n^c) < \varepsilon$. Thus, $\mathcal{I}$ is tight. It follows, by the Prokhorov Theorem, that $\mathcal{I}$ is relatively compact. Combining relative compactness with the previously shown closedness, we have that $\mathcal{I}$ is compact. Hence, by Krein-Milman Theorem, $\mathcal{I}$ has at least one extreme point.
\end{proof}}

Now, we can state the existence of an ergodic invariant measure. 
\begin{theorem}\label{thm:erginv}
        Let\;\Cref{ass:cv} and the inequality 
    \begin{equation}
        \frac{\nu_0\lambda_1}{2} + 2\nu_1\tilde{C} > K
    \end{equation}
    hold. Then, the transition probability function $P_t$ corresponding to\;\eqref{eq:af} admits at least one ergodic invariant measure.
\end{theorem}

\begin{proof}
    By\;\Cref{thm:inv}, we have there exists at least one invariant measure, so that $\mathcal{I}$ is nonempty.\ \Cref{prop:extremal} tells us that there exists at least one extremal point in the space of invariant measures. It follows by\;\Cref{prop:dpz} that this extremal point in $\mathcal{I}$ is ergodic. Hence, there exists at least one ergodic invariant measure. 
\end{proof}

\section{Uniqueness of Ergodic Invariant Measure}
In this section, we establish the uniqueness of the ergodic invariant measure. We do so by following the recipe outlined in\;\cite{ghmr}. Since our noise is multiplicative adhering to\;\Cref{ass:cv}, we also take into account the techniques of Ferrario and Zanella\;\cite{ferr}.

Given $\sigma$ satisfying\;\Cref{ass:cv}, let $u$ denote the solutions of\;\eqref{eq:af} with $u_0\in H$. Given $\lambda > 0$ and $N> 0$, let $v$ denote the solution of the following with $v_0\in H$.
\begin{equation}\label{eq:nudgedaf}
    v(t) + \int_0^t (\nu_0 Av + \nu_1\mathcal{A}_1(v) + B(v))ds = v_0 + \int_0^t f ds + \lambda\int_0^t P_N(u-v) ds + \int_0^t \sigma(v) dW(s)
\end{equation} 
where $P_N$ denotes the orthogonal projection from $H$ onto the space $\text{span} {\{e_n\}}_{1\leq n\leq N}$. The parameters $\lambda$ and $N$ will be chosen later. We will refer to\;\eqref{eq:nudgedaf} as the nudged equation corresponding to equation\;\eqref{eq:af}. The additional term, $\lambda P_N(u-v)$, does not significantly affect the local monotonicity and the a priori estimates of\;\eqref{eq:af}. Hence, the well-posedness of\;\eqref{eq:nudgedaf} follows similar to\;\eqref{eq:af}.

Let us begin with the moment estimates of the nudged equation.
\begin{proposition}[Moment Estimates for Nudged Equation]\label{prop:momentnudged}
    Let $\sigma$ satisfy\;\Cref{ass:cv} and assume 
    \[\frac{\nu_0\lambda_1}{2} + 2\nu_1\tilde{C} > K.\]
    For $p\in \left[2, 1+ \frac{\nu_0\lambda_1 + 4\tilde{C}\nu_1}{2K}\right)$, the solution $v$ of the nudged NSNE with shear-thickening viscosity satisfies 
    \begin{equation}\label{eq:nudgedmoment}
        \mathbb{E}[|v(t)|^p] \leq |v_0|^p e^{-\overline{a}t} + \frac{\overline{C}}{\overline{a}}(1+|u_0|^2) \qquad \forall t>0
    \end{equation}
    where $\overline{a} = \overline{a}(p,\nu_0,\nu_1,\tilde{C},\lambda_1, K)>0$ and $\overline{C} = \overline{C}(p,\nu_0,\nu_1,\tilde{C}, \lambda_1,\lambda, ||f||_{V'}, K)>0$.
\end{proposition}

\begin{proof}
    For $p=2$, we get the equality
    \begin{equation}\label{eq:nudgedeq1}
    \begin{aligned}
        |v(t)|^2 + 2\nu_0\int_0^t ||v(s)||^2 ds &+ 2\nu_1 \int_0^t ||v(s)||_{V_{1,q}}^q ds = |v_0|^2 + 2\int_0^t \langle f, v(s)\rangle ds\\
        & + 2\lambda\int_0^t \langle P_N(u(s)-v(s)),v(s)\rangle ds + \int_0^t ||\sigma(v(s))||_{L_Q}^2 ds + 2\int_0^t \langle v, \sigma(v(s)) dW(s)\rangle.
    \end{aligned}
    \end{equation}
    Observe that 
    \begin{align*}
        2\lambda \langle P_N (u-v),v\rangle &= 2\lambda (\langle P_N(u), v\rangle - |P_N v|^2)\\
        &\leq 2\lambda |u| \ |v| \\
        &\leq 2\lambda \bigg(\frac{\lambda}{\nu_0 \lambda_1} |u|^2 + \frac{\nu_0\lambda_1}{4\lambda}|v|^2\bigg)\\
        &= \frac{2\lambda^2}{\nu_0\lambda_1}|u|^2 + \frac{\nu_0\lambda_1}{2} |v|^2.
    \end{align*}
    Applying this estimate, the Young inequality, and Assumption\;\ref{M2} we get
    \begin{equation*}
    \begin{aligned}
        |v(t)|^2 + \nu_0\int_0^t ||v(s)||^2 ds &+ 2\nu_1\int_0^t ||v(s)||^q_{V_{1,q}} ds \leq |v_0|^2 + \frac{1}{\nu_0}||f||^2_{V'}t+ \frac{2\lambda^2}{\nu_0\lambda_1}\int_0^t |u(s)|^2 ds\\
        & + \frac{\nu_0\lambda_1}{2} \int_0^t |v(s)|^2 ds + Kt + K\int_0^t |v(s)|^2 ds + 2\int_0^t \langle v(s), \sigma(v(s)) dW(s)\rangle.
    \end{aligned}
    \end{equation*}
    Now, applying the Poincar\'{e} inequality and take the expectation yields 
    \begin{align*}
        \mathbb{E}[|v(t)|^2] + \bigg(\frac{\nu_0\lambda_1}{2} + 2\nu_1\tilde{C} - K\bigg)\int_0^t \mathbb{E}[|v(s)|^2] ds \leq |v_0|^2 + \bigg(\frac{1}{\nu_0} ||f||_{V'}^2 &+ 2\nu_1\tilde{C} + K\bigg)t\\
        & + \frac{2\lambda^2}{\nu_0\lambda_1} \int_0^t \mathbb{E}[|u(s)|^2] ds.
    \end{align*}
    By\;\Cref{prop:moment}, we have $\sup_{t\geq 0} \mathbb{E}[|u|^2] \leq C_1(1 + |u_0|^2)$ for a suitable constant $C_1$ which gives us 
    \begin{align*}
        \mathbb{E}[|v(t)|^2] + \bigg(\frac{\nu_0\lambda_1}{2} + 2\nu_1\tilde{C} - K\bigg)\int_0^t \mathbb{E}[|v(s)|^2] ds \leq |v_0|^2 + \bigg(\frac{1}{\nu_0} ||f||_{V'}^2 &+ 2\nu_1\tilde{C} + K\bigg)t\\
        & + \frac{2\lambda^2}{\nu_0\lambda_1}C_1 (1+|u_0|^2)t.
    \end{align*}
    Rewriting this estimate with condensed coefficients, we say 
    \begin{align*}
        \mathbb{E}[|v(t)|^2] + \overline{a}\int_0^t \mathbb{E}[|v(s)|^2] ds \leq |v_0|^2 + \overline{C}(1+|u_0|^2)t
    \end{align*}
    where 
    \[\overline{a} = \frac{\nu_0\lambda_1}{2} + 2\nu_1\tilde{C} - K\]
    and $\overline{C} = \overline{C}(\nu_0,\nu_1, \tilde{C}, \lambda_1,\lambda, ||f||_{V'}, K)$. Now, applying the Gr\"{o}nwall inequality, we have
    \[\mathbb{E}[|v(t)|^2] \leq |v_0|^2 e^{-\overline{a} t} + \frac{\overline{C}}{\overline{a}}(1+|u_0|^2).\]
    Now, let $p\in \left(2, 1+ \frac{\nu_0\lambda_1 + 4\tilde{C}\nu_1}{2K}\right)$. By the It\^{o} formula in differential form, we get 
    \begin{equation}\label{eq:nudgedpdif}
        \begin{aligned}
            d|v(t)|^p + (p\nu_0 |v(t)|^{p-2} ||v(t)||^2 &+ p\nu_1 |v(t)|^{p-2} ||v(t)||_{V_{1,q}}^q)dt \leq p|v(t)|^{p-2}\langle f,v(t)\rangle dt\\
            &\quad + \frac{p(p-1)}{2} |v(t)|^{p-2} ||\sigma(v(t))||_{L_Q}^2 dt + \lambda p|v(t)|^{p-2}\langle P_N(u(t)-v(t)), v(t)\rangle dt\\
            &\quad + p|v(t)|^{p-2}\langle v(t), \sigma(v(t)) dW(t)\rangle.    
        \end{aligned}
    \end{equation}
    Observe the following 
    \begin{align*}
        p|v|^{p-2}\langle f,v \rangle &\leq p|v|^{p-2}||v||\ ||f||_{V'}\\
        &\leq p|v|^{p-2}\bigg( \frac{\nu_0}{2} ||v||^2 + \frac{1}{2\nu_0}||f||^2_{V'}\bigg)\\
        &= \frac{p\nu_0}{2}|v|^{p-2}||v||^2 + \frac{p}{2\nu_0} |v|^{p-2}||f||^2_{V'}\\
        &\leq \frac{p\nu_0}{2}|v|^{p-2}||v||^2 + \frac{p}{2\nu_0}\bigg(\frac{\nu_0^2\lambda_1(p-2)}{4p}|v|^p + \frac{2^{p-1}}{p\nu_0^{p-2}\lambda_1^{\frac{p}{2} - 1}}||f||_{V'}^p\bigg)\\
        &= \frac{p\nu_0}{2}|v|^{p-2}||v||^2 + \frac{\nu_0\lambda_1(p-2)}{8}|v|^p + \frac{2^{p-2}}{\nu_0^{p-1} \lambda_1^{\frac{p}{2}-1}}||f||_{V'}^p
    \end{align*} 
    \begin{align*}
        \frac{p(p-1)}{2} |v|^{p-2} ||\sigma(v)||_{L_Q}^2 &\leq \frac{p(p-1)}{2}|v|^{p-2}K(1+|v|^2)\\
        &= \frac{p(p-1)}{2}|v|^{p-2}K + \frac{Kp(p-1)}{2}|v|^{p}\\
        &\leq \frac{p(p-1)}{2}\bigg(\frac{\nu_0\lambda_1}{2p(p-1)}|v|^p + \frac{2^{\frac{p}{2}} }{p}{\Big(\frac{(p-1)(p-2)}{\nu_0\lambda_1}\Big)}^{\frac{p}{2}-1}K^{\frac{p}{2}}\bigg)\\
        &\qquad + \frac{Kp(p-1)}{2}|v|^{p}\\
        &= \bigg(\frac{\nu_0\lambda_1}{4} + \frac{Kp(p-1)}{2}\bigg)|v|^p + \frac{2^{\frac{p}{2}-1}{(p-1)}^{\frac{p}{2}}{(p-2)}^{\frac{p}{2}-1}}{{(\nu_0\lambda_1)}^{\frac{p}{2}-1}}K^{\frac{p}{2}}
    \end{align*}
    and
    \begin{align*}
        \lambda p |v|^{p-2} \langle P_N(u-v),v\rangle &\leq \lambda p |v|^{p-2}|u|\ |v|\\
        &\leq \lambda p |v|^{p-2}\bigg(\frac{4\lambda}{\nu_0\lambda_1}|u|^2 + \frac{\nu_0\lambda_1}{16\lambda}|v|^2\bigg)\\
        &= \frac{4\lambda^2 p}{\nu_0\lambda_1} |v|^{p-2}|u|^2 + \frac{p\nu_0\lambda_1}{16}|v|^p\\
        &\leq \frac{4\lambda^2 p}{\nu_0\lambda_1}\bigg({\Big(\frac{\nu_0\lambda_1}{8\lambda}\Big)}^2 |v|^p + \frac{2{(p-2)}^{\frac{p}{2}-1}{(8\lambda)}^{p-2}}{p^{\frac{p}{2}}{(\nu_0\lambda_1)}^{p-2}}|u|^p\bigg) + \frac{p\nu_0\lambda_1}{16}|v|^p\\
        &= \frac{p\nu_0\lambda_1}{8}|v|^p + \frac{8^{p-1}\lambda^p{(p-2)}^{\frac{p}{2}-1}}{p^{\frac{p}{2}-1} {(\nu_0\lambda_1)}^{p-1}}|u|^p.
    \end{align*}
    Applying these inequalities and\;\Cref{prop:lowerbound} to\;\eqref{eq:nudgedpdif}, we get 
    \begin{align*}
        d|v(t)|^p + &(p\nu_0 |v|^{p-2} ||v(t)||^2 + \tilde{C}p\nu_1(|v(t)|^p - 1))dt\\
        &\leq \bigg(\frac{p\nu_0}{2}|v(t)|^{p-2}||v(t)||^2 + \bigg(\frac{p\nu_0\lambda_1}{4} + \frac{Kp(p-1)}{2}\bigg)|v(t)|^p + \frac{2^{p-2}}{\nu_0^{p-1} \lambda_1^{\frac{p}{2}-1}}||f||_{V'}^p \\
        &+ \frac{2^{\frac{p}{2}-1}{(p-1)}^{\frac{p}{2}}{(p-2)}^{\frac{p}{2}-1}}{{(\nu_0\lambda_1)}^{\frac{p}{2}-1}}K^{\frac{p}{2}} + \frac{8^{p-1}\lambda^p{(p-2)}^{\frac{p}{2}-1}}{p^{\frac{p}{2}-1} {(\nu_0\lambda_1)}^{p-1}}|u(t)|^p\bigg)dt + p|v(t)|^{p-2}\langle v(t), \sigma(v(t)) dW(t)\rangle.
    \end{align*}
    Applying the Poincar\'{e} inequality and simplifying gives us 
    \begin{align*}
        d|v(t)|^p + \bigg(\frac{p\nu_0\lambda_1}{4} +\tilde{C}p\nu_1 &- \frac{Kp(p-1)}{2} \bigg) |v(t)|^p dt \leq \bigg(\tilde{C}p\nu_1 + \frac{2^{p-2}}{\nu_0^{p-1} \lambda_1^{\frac{p}{2}-1}}||f||_{V'}^p\\
        &+ \frac{2^{\frac{p}{2}-1}{(p-1)}^{\frac{p}{2}}{(p-2)}^{\frac{p}{2}-1}}{{(\nu_0\lambda_1)}^{\frac{p}{2}-1}}K^{\frac{p}{2}}+ \frac{8^{p-1}\lambda^p{(p-2)}^{\frac{p}{2}-1}}{p^{\frac{p}{2}-1} {(\nu_0\lambda_1)}^{p-1}}|u|^p\bigg)dt\\
        &+ p|v(t)|^{p-2}\langle v(t), \sigma(v(t)) dW(t)\rangle.
    \end{align*}
    Taking the expectation yields
    \[\frac{d}{dt} \mathbb{E}[|v(t)|^p] \leq -\overline{a} \mathbb{E}[|v(t)|^p] + C_2 + C_3\mathbb{E}[|u(t)|^p]\]
    where 
    \begin{align*}
        \overline{a} &= \frac{p\nu_0\lambda_1}{4} +\tilde{C}p\nu_1 - \frac{Kp(p-1)}{2}\\
        C_2 &= \tilde{C}p\nu_1 + \frac{2^{p-2}}{\nu_0^{p-1} \lambda_1^{\frac{p}{2}-1}}||f||_{V'}^p + \frac{2^{\frac{p}{2}-1}{(p-1)}^{\frac{p}{2}}{(p-2)}^{\frac{p}{2}-1}}{{(\nu_0\lambda_1)}^{\frac{p}{2}-1}}K^{\frac{p}{2}}\\
        C_3 &= \frac{8^{p-1}\lambda^p{(p-2)}^{\frac{p}{2}-1}}{p^{\frac{p}{2}-1} {(\nu_0\lambda_1)}^{p-1}}.
    \end{align*}
    By\;\Cref{prop:moment}, we have $\sup_{t\geq 0}\mathbb{E}[|u(t)|^p] \leq C_4(1+|u_0|^p)$ for a suitable constant $C_4$. It follows that 
    \[\frac{d}{dt} \mathbb{E}[|v(t)|^p] \leq -\overline{a} \mathbb{E}[|v(t)|^p] + \overline{C}(1 +|u_0|^p)\]
    where $\overline{C} = \overline{C}(p, \nu_0,\nu_1, \tilde{C}, \lambda_1,\lambda, ||f||_{V'}, K)$. Applying the Gr\"{o}nwall inequality, we arrive at 
    \[\mathbb{E}[|v(t)|^p] \leq |v_0|e^{-\overline{a}t} + \frac{\overline{C}}{\overline{a}}(1 + |u_0|^p).\]
\end{proof}
Now, we examine the difference between the solutions of\;\eqref{eq:af} and\;\eqref{eq:nudgedaf}.
\begin{proposition}[Difference Estimate]\label{prop:difference}
    Let $\sigma$ satisfy\;\Cref{ass:cv} and assume 
    \[\frac{\nu_0\lambda_1}{2} + 2\nu_1\tilde{C} > K.\]
    Let $u$ and $v$ denote the solutions to\;\eqref{eq:af} and\;\eqref{eq:nudgedaf} respectively. If we let $\lambda = \frac{\nu_0\lambda_N}{2}$ in\;\eqref{eq:nudgedaf}, then, for any initial conditions $u_0, v_0 \in H$ respectively, we have 
    \begin{equation}\label{eq:diffN}
        \mathbb{E}\Big[e^{\Gamma(t)} |u(t) - v(t)|^2\Big] \leq |u_0 - v_0|^2
    \end{equation}
    where 
    \begin{equation}\label{eq:gammaN}
      \Gamma(t) = (\nu_0\lambda_N - L)t - \frac{2}{\nu_0}\int_0^t ||u(s)||^2 ds.  
    \end{equation}
\end{proposition}

\begin{proof}
    Let $w(t) := u(t)-v(t)$ and $w_0 := u_0 - v_0$. Applying the It\^{o} formula to $|w(t)|^2$, we get
    \begin{align*}
        d|w(t)|^2 + (2\nu_0||w(t)||^2 + 2\nu_1\langle \mathcal{A}_1(u(t)) &- \mathcal{A}_1(v(t)), w(t)\rangle + 2\langle B(u(t)) - B(v(t)), w(t)\rangle + 2\lambda |P_N w(t)|^2)dt\\
        &=||\sigma(u(t))-\sigma(v(t))||^2_{L_Q}dt + 2\langle w(t), (\sigma(u(t)) - \sigma(v(t)))dW(t) \rangle
    \end{align*}
    Observe the following inequalities 
    \begin{align*}
        \langle \mathcal{A}_1(u) - \mathcal{A}_1(v), w\rangle &\geq \gamma ||w||^q_{V_{1,q}}\\
        |\langle B(u) - B(v), w \rangle| &\leq \frac{\nu_0}{2}||w||^2 + \frac{1}{\nu_0} ||u||^2 |w|^2.
    \end{align*}
    Applying these inequalities and Assumption\;\ref{M3} gives us 
    \begin{align*}
        d|w(t)|^2 + (\nu_0 ||w(t)||^2 + 2\nu_1\gamma||w(t)||^q_{V_{1,q}} + 2\lambda |P_N w(t)|^2)dt \leq \bigg(\frac{2}{\nu_0}||u(t)||^2+ L\bigg)|w(t)|^2 dt + dM(t)
    \end{align*}
    where $M(t) = 2\int_0^t \langle w(s), (\sigma(u(s)) - \sigma(v(s)))dW(s) \rangle$. For the control term, we have, by the generalized Poincar\'{e} inequality,
    \begin{align*}
        \nu_0||w||^2 + 2\lambda |P_N w|^2 &\geq \nu_0||Q_N w||^2 + 2\lambda |P_N w|^2\\
        &\geq \nu_0\lambda_N|Q_N w|^2 + 2\lambda|P_N w|^2.
    \end{align*}
    Let $\lambda = \frac{\nu_0\lambda_N}{2}$. Then 
    \begin{align*}
        \nu_0||w||^2 + 2\lambda |P_N w|^2 &\geq \nu_0\lambda_N|Q_N w|^2 + \nu_0\lambda_N |P_N w|^2\\
        &= \nu_0\lambda_N |w|^2.
    \end{align*}
    Substituting this and removing $2\nu_1\gamma||w||^q_{V_{1,q}}$ (since it is positive), we have 
    \begin{align*}
        d|w(t)|^2 + \bigg(\nu_0 \lambda_N - \frac{2}{\nu_0}||u(t)||^2 - L \bigg)|w(t)|^2dt \leq dM(t).
    \end{align*}
    Let 
    \begin{equation*}
        \Gamma(t) = (\nu_0\lambda_N - L)t - \frac{2}{\nu_0}\int_0^t ||u(s)||^2 ds.
    \end{equation*}
    Then multiplying by $e^{\Gamma(t)}$ and integrating gives
    \begin{align*}
        e^{\Gamma(t)}|w(t)|^2 - |w_0|^2 \leq \int_0^t e^{\Gamma(s)}dM(s).
    \end{align*}
    Taking expectation yields 
    \[\mathbb{E}\big(e^{\Gamma(t)} |w(t)|^2\big) \leq |w_0|^2\]
    as desired.
\end{proof}

Consider the following
\begin{equation}\label{eq:stoppingtime}
    \tau_{\beta, R} := \inf\bigg\{t\geq 0 ; \frac{2}{\nu_0} \int_0^t ||u(s)||^2 ds + \Big(L - \frac{\nu_0\lambda_N}{2}\Big)t - \beta \geq R\bigg\}
\end{equation}
for $\beta, R > 0$.

\begin{remark}\label{remark:GammaN}
    On the event $\{\tau_{\beta,R} = \infty\}$, we have 
    \[\frac{2}{\nu_0} \int_0^t ||u(s)||^2 ds + \Big(L -\frac{\nu_0 \lambda_N}{2}\Big)t - \beta < R\]
    for $t\geq 0$. It follows, for $\Gamma(t)$ as defined in\;\eqref{eq:gammaN}, that 
    \[\frac{\nu_0\lambda_N}{2}t - \beta - R < \Gamma(t)\]
    on the event $\{\tau_{\beta,R} = \infty\}$ for all $t\geq 0$.
\end{remark}

\begin{proposition}[Estimate in Probability]\label{prop:estinprob}
    Assume\;\Cref{ass:cv} holds. Let $u$ denote the solution of the NSNE with shear-thickening viscosity\;\eqref{eq:af}. If 
    \begin{equation}\label{eq:extra}
        \frac{\nu_0\lambda_1}{4} + \tilde{C}\nu_1 > K,
    \end{equation}
    then, for any $p \in (2, \frac{1}{2} + \frac{\nu_0\lambda_1 + 4\tilde{C}\nu_1}{4k})$, 
    \[\mathbb{P}\bigg(\sup_{t\geq 0} \bigg[|u(t)|^2 + a\int_0^t ||u(s)||^2 ds + 2\nu_1\int_0^t ||u(s)||_{V_{1,q}}^q ds - |u_0|^2 - (b + 1)(t+1)\bigg] \geq R\bigg) \leq \frac{C(1 + {|u_0|}^{2p})}{R^{\frac{p}{2} - 1}}\]
    for all $R > 0$ where $a = \nu_0 - \frac{K}{\lambda_1}$, $b = \frac{1}{\nu_0}||f||_{V'}^2 + K$, and $C = C(p, \nu_0, \nu_1, \tilde{C}, \lambda_1, ||f||_{V'}, K) > 0$.
\end{proposition}

\begin{proof}
    Using similar techniques as those used in the moment estimates, we get 
    \begin{align*}
        |u(t)|^2 + a\int_0^t ||u(s)||^2 ds + 2\nu_1\int_0^t ||u(s)||_{V_{1,q}}^q ds - |u_0|^2 - bt \leq M(t)
    \end{align*} 
    where $a = \nu_0 - \frac{K}{\lambda_1}$, $b = \frac{1}{\nu_0}||f||_{V'}^2 + K$, and $M(t) = 2\int_0^t \langle u(s), \sigma(u(s)) dW(s)\rangle$. Subtracting $b+1 + t$ on both sides gives us 
    \begin{align*}
        |u(t)|^2 + a\int_0^t ||u(s)||^2 ds + 2\nu_1\int_0^t ||u(s)||_{V_{1,q}}^q ds - |u_0|^2 - (b+1)(t+1) \leq M(t) - t - (b+1).
    \end{align*}
    It follows that 
    \begin{align*}
        \mathbb{P}\bigg(\sup_{t\geq 0} \bigg[|u(t)|^2 + a\int_0^t ||u(s)||^2 ds + 2\nu_1\int_0^t ||u(s)||_{V_{1,q}}^q ds &- |u_0|^2 - (b + 1)(t+1)\bigg] \geq R\bigg)\\
        &\leq \mathbb{P}\Big(\sup_{t\geq 0} [M(t) - t - (b+1)] \geq R\Big).
    \end{align*}
    Observe that 
    \[\Big\{\sup_{t\geq 0} [M(t) - t - (b+1)] \geq R\Big\} = \bigcup_{m\geq 0} \Big\{\sup_{t\in [m,m+1)} [M(t) - t - (b+1)] \geq R\Big\}.\]
    We then say 
    \[\Big\{\sup_{t\in [m,m+1)} [M(t) - t - (b+1)] \geq R\Big\} \subset \Big\{\sup_{t\in [0, m+1]} |M(t)| \geq R + m + b+1\Big\}.\]
    So, it follows that 
    \begin{equation}\label{eq:estprob}
    \begin{aligned}[b]
        \mathbb{P}\bigg(\sup_{t\geq 0} \bigg[|u(t)|^2 + a\int_0^t ||u(s)||^2 ds &+ 2\nu_1\int_0^t ||u(s)||_{V_{1,q}}^q ds - |u_0|^2 - (b + 1)(t+1)\bigg] \geq R\bigg)\\
        &\leq \sum_{m\geq 0} \mathbb{P}\bigg(\sup_{t\in [0, m+1]} |M(t)| \geq R + m + b+1\bigg)\\
        &\leq \sum_{m\geq 0} \frac{\mathbb{E}\Big({(\sup_{t\in [0,m+1]}|M(t)|)}^p\Big)}{{(R+m+b+1)}^p}
    \end{aligned}
    \end{equation}
    where the last inequality is by the Markov inequality. To continue the estimate on this probability, we must, first, estimate the expectation. We will use the Burkholder-Davis-Gundy inequality to do this. Observe that
    \begin{equation}\label{eq:bdg2}
    \begin{aligned}[b]
        \mathbb{E}\Big({(\sup_{t\in [0,m+1]}|M(t)|)}^p\Big) &\leq C\mathbb{E}\Big( {({[M]}_{m+1})}^{\frac{p}{2}}\Big)\\
        &\leq C\mathbb{E}\bigg( {\Big(\int_0^{m+1} (1+|u(s)|^4) ds\Big)}^{\frac{p}{2}}\bigg)\\
        &\leq C {(m+1)}^{\frac{p}{2} - 1}\mathbb{E}\bigg(\int_0^{m+1} (1 + |u(s)|^{2p}) ds\bigg)\\
        &\leq C {(m+2)}^{\frac{p}{2}} (1+ |u_0|^{2p})
    \end{aligned}
    \end{equation}
    where the first step is by the Burkholder-Davis-Gundy inequality, the second by Assumption\;\ref{M2}, third by the H\"{o}lder inequality, and the last by\;\Cref{prop:moment}. Note that, to use\;\Cref{prop:moment}, we must have $2p < 1 + \frac{\nu_0\lambda_1 + 4\tilde{C}\nu_1}{2K}$. Also, note that the constant $C$ may be different from step to step. Applying\;\eqref{eq:bdg2} to\;\eqref{eq:estprob}, we get 
    \begin{align*}
        \mathbb{P}\bigg(\sup_{t\geq 0} \bigg[|u(t)|^2 + a\int_0^t ||u(s)||^2 ds &+ 2\nu_1\int_0^t ||u(s)||_{V_{1,q}}^q ds - |u_0|^2 - (b + 1)(t+1)\bigg] \geq R\bigg)\\
        &\leq C(1+ |u_0|^{2p})\sum_{m\geq 0} \frac{{(m+2)}^{\frac{p}{2}}}{{(R+m+b+1)}^p}.
    \end{align*}
    This series converges precisely when $p > 2$ and we arrive at  
    \[\mathbb{P}\bigg(\sup_{t\geq 0} \bigg[|u(t)|^2 + a\int_0^t ||u(s)||^2 ds + 2\nu_1\int_0^t ||u(s)||_{V_{1,q}}^q ds - |u_0|^2 - (b + 1)(t+1)\bigg] \geq R\bigg)\leq \frac{C(1+ {|u_0|}^{2p})}{R^{\frac{p}{2} -1}}.\]
    To be able to use\;\Cref{prop:moment} and have the series converge, we require that $p\in (2,  \frac{1}{2} + \frac{\nu_0\lambda_1 + 4\tilde{C}\nu_1}{4K})$ which in turn requires $2 < \frac{1}{2} + \frac{\nu_0\lambda_1 + 4\tilde{C}\nu_1}{4K}$, but this is precisely\;\eqref{eq:extra}. Thus, our conclusion has been proven.
\end{proof}

\begin{proposition}[Decay Estimate]\label{prop:decay}
    Assume\;\Cref{ass:cv} and\;\eqref{eq:extra} hold. Consider the stopping time $\tau_{\beta,R}$ defined in\;\eqref{eq:stoppingtime}, where $u$ is the solution to the NSNE with shear-thickening viscosity\;\eqref{eq:af}. If 
    \[\beta \geq \frac{2(b+1+|u_0|^2)}{\nu_0(\nu_0 - \frac{K}{\lambda_1})}\]
    with $b = \frac{1}{\nu_0}||f||_{V'}^2 + K$, then there exists an integer $\overline{N} > 0$ such that for any $N \geq \overline{N}$,
    \[\mathbb{P}(\tau_{\beta,R} < \infty) \leq \frac{C(1+|u_0|^{4(p+1)})}{R^p}\]
    for any $p \in (0, \frac{\nu_0\lambda_1 + 4\tilde{C} \nu_1}{8K} - \frac{3}{4})$ where $C > 0$ is a constant independent of $\beta, R$, and $u_0$.  
\end{proposition}

\begin{proof}
    Define the set 
    \[A_{\beta, R} = \bigg\{\sup_{t\geq 0} \bigg(\frac{2}{\nu_0} \int_0^t ||u(s)||^2 ds + \Big(L - \frac{\nu_0\lambda_N}{2}\Big)t - \beta\bigg)\geq R\bigg\}.\]
    It is clear that $\mathbb{P}(\tau_{\beta,R} < \infty) \leq \mathbb{P}(A_{\beta,R})$. So, we need to estimate $\mathbb{P}(A_{\beta,R})$. Observe that the complement of $A_{\beta,R}$ can be written as 
    \[A_{\beta,R}^C = \bigg\{ \Big(\nu_0 - \frac{K}{\lambda_1}\Big)\int_0^t ||u(s)||^2 ds < \frac{\nu_0}{2}\Big(\nu_0- \frac{K}{\lambda_1}\Big)\Big[R + \beta +\Big(\frac{\nu_0\lambda_N}{2} - L\Big)t\Big] \quad \text{ for any } t\geq 0\bigg\}.\]
    Take $\overline{N}$ large enough such that for any $N\geq \overline{N}$
    \[\lambda_N > \frac{2L}{\nu_0} + \frac{4(b+1)}{\nu_0^2(\nu_0 - \frac{K}{\lambda_1})},\]
    choose 
    \[\beta\geq \frac{2(b+1+|u_0|^2)}{\nu_0(\nu_0 - \frac{K}{\lambda_1})},\]
    and set $\overline{R} := \nu_0(\nu_0 - \frac{K}{\lambda_1})R$. Then 
    \[A_{\beta,R}^C \supseteq \bigg\{\Big(\nu_0 - \frac{K}{\lambda_1}\Big)\int_0^t ||u(s)||^2 ds < \overline{R} + |u_0|^2 + (b+1)(t+1) \quad \text{ for any } t\geq 0\bigg\}.\]
    It follows that 
    \[A_{\beta, R} \subseteq \bigg\{\sup_{t\geq 0} \bigg[\Big(\nu_0 - \frac{K}{\lambda_1}\Big)\int_0^t ||u(s)||^2 ds - |u_0|^2 - (b+1)(t+1)\bigg] \geq \overline{R}\bigg\}.\]
    By\;\Cref{prop:estinprob}, we have  
    \[\mathbb{P}(A_{\beta,R}) \leq \frac{C(1+|u_0|^{2p})}{\overline{R}^{\frac{p}{2} - 1}}\]
    provided $\frac{\nu_0\lambda_1}{4} + \tilde{C}\nu_1 > K$ and $p \in (2, \frac{1}{2} + \frac{\nu_0\lambda_1 + 4\tilde{C}\nu_1}{4K})$ where $C>0$ is independent of $\beta, R$, and $u_0$. Allowing the constant $C$ to absorb the constant in the definition of $\overline{R}$ and rescaling our exponents, we say 
    \[\mathbb{P}(\tau_{\beta,R} < \infty) \leq \frac{C(1+|u_0|^{4(p+1)})}{R^p}\]
    where $p \in (0, \frac{\nu_0\lambda_1 + 4\tilde{C} \nu_1}{8K} - \frac{3}{4})$.
\end{proof}

To prove uniqueness of invariant measure, we introduce a modified version of the\;\eqref{eq:af} to form a generalized coupling. We want this modified equation to be equivalent in law to equation\;\eqref{eq:af} and, for any distinct initial conditions, there is a positive probability that the solutions converge as time goes to infinity. Before defining such an equation, we will fix two initial conditions $u_0, v_0 \in H$ of the equations\;\eqref{eq:af} and\;\eqref{eq:nudgedaf} respectively. Define the oppositional control term $h$ by
\begin{equation}\label{eq:h}
    h(t) := \lambda g(v(t)) P_N(u(t) - v(t)), \quad t\geq 0
\end{equation}
where $\lambda P_N(u(t)-v(t))$ is the nudged term in equation\;\eqref{eq:nudgedaf} and $g$ respects the following assumptions.
\begin{assumption}\label{ass:g}
  There exists a measurable map $g: H\rightarrow L(H,U)$ such that 
  \begin{equation}\label{eq:g1}
    \sup_{u\in H}||g(u)||_{L(H,U)} < \infty
  \end{equation}
  and 
  \begin{equation}
    \sigma(u)g(u) = P_M,
  \end{equation}
  for all $u\in H$, for a positive integer $M$.
\end{assumption}
By these assumptions, the definition\;\eqref{eq:h} of $h$ makes sense and we can see that $h$ lives in the space $U$ with the noise. Next, we define the stopping time $\theta_\kappa$, for a given $\kappa > 0$, by 
\begin{equation}\label{eq:theta}
    \theta_\kappa := \inf\bigg\{t\geq 0 : \int_0^t |P_N(u(s) - v(s))|^2 ds  \geq \kappa\bigg\}.
\end{equation} 
As in definition\;\eqref{eq:nudgedaf} of $v$, the $N$ in\;\eqref{eq:h} and\;\eqref{eq:theta} is fixed. The constant $\kappa$ is also fixed. These fixed constants will be determined later. The modified equation we will use to form a generalized coupling is given by 
\begin{equation}\label{eq:stoppednudged}
    \tilde{v}(t) + \int_0^t (\nu_0 A\tilde{v} + \nu_1\mathcal{A}_1(\tilde{v}) + B(\tilde{v})) ds = v_0 + \int_0^t f ds + \int_0^t \sigma(\tilde{v}) d\widetilde{W}(s)  
\end{equation}
Here $\widetilde{W}$ is given by 
\[\widetilde{W}(t) := W(t) + \int_0^t h(s)\mathbbm{1}_{\{s\leq \theta_\kappa\}} ds.\]

Equation\;\eqref{eq:stoppednudged} will be referred to as the stopped nudged equation corresponding to equation\;\eqref{eq:af}. Note that on the set $\{\theta_\kappa = \infty\}$, $\tilde{v} = v $ $\mathbb{P}$-a.s., where $v$ denotes the solution to the nudged equation\;\eqref{eq:nudgedaf}.

\begin{lemma}\label{lemma:girsanov}
    Let\;\Cref{ass:cv} hold. If\;\Cref{ass:g} holds with $M\geq N$, then, for any $\kappa, \lambda > 0$, the laws of solution to\;\eqref{eq:af} and\;\eqref{eq:stoppednudged} are equivalent, i.e. $\mathbb{P}_{v_0} \sim \tilde{\mathbb{P}}_{v_0}$ as measures on $C([0,\infty); H)$.
\end{lemma}

\begin{proof}
    By the definition of $h$ and the H\"{o}lder inequality, we have 
    \begin{align*}
        \int_0^\infty ||h(s)||_{U}^2\, \mathbbm{1}_{s\leq \theta_\kappa} ds &\leq \lambda^2 \Big(\sup_{u\in H} ||g(u)||_{L(H,U)}^2\Big)\int_0^\infty |P_N(u(s)-v(s))|^2\, \mathbbm{1}_{\{s\leq\theta_\kappa\}} ds\\
        &\leq \lambda^2 \Big(\sup_{u\in H} ||g(u)||_{L(H,U)}^2\Big)\kappa\\
        &< \infty.
    \end{align*}
    It follows that $h(s)\,\mathbbm{1}_{\{s\leq\theta_\kappa\}}$ satisfies the Novikov condition 
    \[\mathbb{E}\bigg[\text{exp}\Big(\frac{1}{2}\int_0^\infty ||h(s)||_{U}^2\ \mathbbm{1}_{\{s\leq \theta_\kappa\}} ds \Big)\bigg] < \infty.\]
    Therefore, by the Girsanov Theorem, there exists a probability measure $\mathbb{Q}$ on $C([0,T); U)$ such that $\widetilde{W}$ is an $U$-valued Wiener process on the time interval $[0,\infty)$ with respect to $\mathbb{Q}$. It follows that the laws of the solutions to\;\eqref{eq:stoppednudged} and\;\eqref{eq:af} with initial condition are equivalent, i.e. $\mathbb{P}_{v_0} \sim \tilde{\mathbb{P}}_{v_0}$ as measures on $C([0,\infty); H)$.
\end{proof}

\begin{proposition}
    Let Assumptions\;\ref{ass:cv} and\;\ref{ass:g} and inequality\;\eqref{eq:extra} hold. Let $u$ be the solution of equation\;\eqref{eq:af} and $\tilde{v}$ the solution of the stopped nudged equation\;\eqref{eq:stoppednudged} with initial conditions $u_0, v_0 \in H$ respectively. Then there exists a positive integer $\overline{N}$ and $\lambda>0$, depending on $\overline{N}$ and $\nu_0$, such that 
    \begin{equation}
        \mathbb{P}\Big(\lim_{n\rightarrow \infty} |u(n) - \tilde{v}(n)|^2 = 0\Big) > 0
    \end{equation}
    when $N\geq \overline{N}$.
\end{proposition}

\begin{proof}
    We begin by introducing the event
    \[B_n := \bigg\{|u(n) - v(n)|^2 + \int_n^{n+1} |P_N (u(s) - v(s))| ds > \frac{1}{n^2}\bigg\}\]
    for any $n > 0$. We set 
    \[B := \bigcap_{m=1}^\infty \bigcup_{n=m}^\infty B_n.\]
    Now, we fix $\overline{N}$ large enough, so that the previous propositions hold. Then, we fix $N \geq \overline{N}$ and set $\lambda = \frac{\nu_0\lambda_N}{2}$ in\;\eqref{eq:nudgedaf}. Observe that the probability of $B$ can be decomposed into 
    \[\mathbb{P}(B) = \mathbb{P}(B \cap \{\tau_{\beta,R} = \infty\}) + \mathbb{P}(B \cap \{\tau_{\beta,R} < \infty\}).\]
    For the first probability, we see, by the Markov inequality and the Fubini Theorem, that 
    \begin{align*}
        \mathbb{P}(B_n \cap \{\tau_{\beta,R} = \infty\}) &\leq n^2 \mathbb{E}\bigg[1_{\{\tau_{\beta,R} =\infty\}} \Big(|u(n) - v(n)|^2 + \int_n^{n+1} |u(s) - v(s)|^2 ds\Big)\bigg]\\
        &\leq n^2\mathbb{E}|w_0|^2 e^{\beta + R - \frac{\nu_0\lambda_N}{2}n} + \frac{2}{\nu_0\lambda_N} n^2\mathbb{E}|w_0|^2 e^{\beta+R -\frac{\nu_0\lambda_N}{2}n}
    \end{align*}
    where the last inequality follows by\;\Cref{prop:difference} and\;\Cref{remark:GammaN}. It follows that $\sum_{n=1}^\infty \mathbb{P}(B_n \cap \{\tau_{\beta,R} = \infty\}) < \infty$. Hence, by the Borel-Cantelli lemma, 
    \[\mathbb{P}(B \cap \{\tau_{\beta,R} = \infty\}) = 0.\]
    Now, for the second probability, we have, for $p \in (0, \frac{\nu_0\lambda_1 + 4\tilde{C} \nu_1}{8K} - \frac{3}{4})$,
    \[\mathbb{P}(B) = \mathbb{P}(B\cap \{\tau_{\beta,R} < \infty\}) \leq \mathbb{P}(\tau_{\beta, R} < \infty) \leq \frac{C}{R^p}\]
    where the last inequality is due to\;\Cref{prop:decay} and $C$ does not depend on $R$. Choosing $R^*$ large enough, we have $\mathbb{P}(B)$ is close to $0$. Hence $\mathbb{P}(B^c)$ is close to $1$. By continuity from below, we can find an $m^* > 0$ large enough so that 
    \[\mathbb{P}\Big(\bigcap_{n = m^*}^\infty B_n^c\Big) > \frac{3}{4}.\]
    For this fixed $m^*$, consider the set 
    \[E_{m^*, R^*} := \bigg\{\int_0^{m^*} |P_N(u(s)-v(s))|^2 ds > R^*\bigg\}.\]
    By the Markov inequality,\;\eqref{eq:moment}, and\;\eqref{eq:nudgedmoment} with $p=2$, we get 
    \begin{align*}
        \mathbb{P}(E_{m^*, R^*}) &\leq \frac{\mathbb{E}\Big(\int_0^{m^*} |P_N(u(s)-v(s))| ds \Big)}{R^*}\\
        &\leq \frac{Cm^*}{R^*}
    \end{align*}
    for some constant $C$ depending on the initial conditions in addition to the parameters mentioned in the propositions. Taking $R^*$ possibly larger, it follows that $\mathbb{P}(E_{m^*,R^*}^c) > \frac{3}{4}$. Hence
    \[\mathbb{P}\Big(E_{m^*,R^*}^c \cap \bigcap_{n=m^*}^\infty B_n^c\Big) > \frac{1}{2}.\]
    Notice on the set $E_{m^*,R^*}^c \cap \bigcap_{n=m^*}^\infty B_n^c$, we have 
    \[\int_0^\infty |P_N(u(s)-v(s))|^2 ds = \int_0^{m^*} |P_N(u(s)-v(s))|^2 ds + \sum_{k=0}^\infty\int_{m^*+k}^{m^*+k+1} |P_N(u(s)-v(s))|^2 ds \leq R^* + \sum_{n=m^*}^\infty \frac{1}{n^2}.\]
    Choose 
    \[\kappa = R^* + \sum_{n=m^*}^\infty \frac{1}{n^2}.\]
    Then $E_{m^*,R^*}^c \cap \bigcap_{n=m^*}^\infty B_n^c \subseteq \{\theta_\kappa = \infty\}$ and, for any $m^*$, $\bigcap_{n=m^*}^\infty B_n^c \subseteq \{\Lim{n\rightarrow \infty} |u(n) - v(n)|^2 = 0\}$. Thus,
    \begin{align*}
        \mathbb{P}\Big(\lim_{n\rightarrow\infty} |u(n) - \tilde{v}(n)|^2 = 0\Big) &\geq \mathbb{P}\bigg(\Big\{\lim_{n\rightarrow\infty} |u(n) - v(n)|^2 = 0\Big\} \cap \{\theta_\kappa = \infty\}\bigg)\\
        &\geq \mathbb{P}\Big(E_{m^*,R^*}^c \cap \bigcap_{n=m^*}^\infty B_n^c\Big)\\
        &> \frac{1}{2}.
    \end{align*}
\end{proof}

\begin{corollary}\label{corollary:posprob}
    Let Assumptions\;\ref{ass:cv} and\;\ref{ass:g} and inequality\:\eqref{eq:extra} hold. Let $u$ be the solution of the equation\;\eqref{eq:af} and $\tilde{v}$ the solution of the stopped nudged equation\;\eqref{eq:stoppednudged} with initial conditions $u_0, v_0 \in H$ respectively. Then there exists a positive integer $\overline{N}$ and $\lambda>0$, depending on $\overline{N}$ and $\nu_0$, such that 
    \begin{equation}\label{eq:stoppeddiffprop}
        \mathbb{P}\Big(\lim_{n\rightarrow \infty} |u(n) - \tilde{v}(n)| = 0\Big) > 0
    \end{equation}
    when $N\geq \overline{N}$.
\end{corollary}

\begin{proof}
    This result follows from the fact that 
    \[\Big\{\lim_{n\rightarrow \infty} |u(n) - \tilde{v}(n)|^2 = 0\Big\} \subset \Big\{\lim_{n\rightarrow \infty} |u(n) - \tilde{v}(n)| = 0\Big\}.\]
\end{proof}

Consider the class of test functions
\[\mathcal{G} := \bigg\{\phi\in C_b(H) : \sup_{u\neq v} \frac{|\phi(u)- \phi(v)|}{|u-v|} < \infty\bigg\}.\]
The set $\mathcal{G}$ is a determining measure set in $H$, i.e., if $\mu_1, \mu_2 \in \mathcal{P}(H)$ are such that $\int_H \phi(u) \mu_1(du) = \int_H \phi(u) \mu_2(du)$ for all $\phi \in \mathcal{G}$, then $\mu_1=\mu_2$. Now, we give the precise statement and proof of our main result.

\begin{theorem}
    Let\;\Cref{ass:cv} and inequality\;\eqref{eq:extra} hold. Then there exists a positive integer $\overline{N}$ such that if\;\Cref{ass:g} holds for some $M\geq \overline{N}$, then $P$ possesses a unique ergodic invariant measure. 
\end{theorem}

\begin{proof}
    Assume there are two ergodic invariant probability measures $\mu_1$ and $\mu_2$ on $H$. Fix any $\phi \in \mathcal{G}$. By the Birkhoff Ergodic Theorem, there exists sets $A_{\mu_1}, A_{\mu_2}\subset H^\mathbb{N}$ such that $\mu_1^\mathbb{N}(A_{\mu_1}) = \mu_2^\mathbb{N}(A_{\mu_2}) = 1$ and if $u \in A_{\mu_1}$ and $v\in A_{\mu_2}$ then 
    \[\lim_{n\rightarrow\infty} \frac{1}{n} \sum_{k=1}^n \phi(u(k)) = \int_H \phi\ d\mu_1 \qquad \text{and} \qquad \lim_{n\rightarrow\infty} \frac{1}{n} \sum_{k=1}^n \phi(v(k)) = \int_H \phi\ d\mu_2.\]
    Define $A_{\mu_1,u_0} := \{\tilde{u} \in A_{\mu_1} : \tilde{u}_0 = u_0\}$, i.e.\ the set of $\tilde{u} \in A_{\mu_1}$ with initial condition $u_0$. Define $A_{\mu_2,v_0}$ analogously. Note that $\mathbb{P}_{u_0}^\mathbb{N}(A_{\mu_1,u_0}) = \mathbb{P}_{u_0}^\mathbb{N}(A_{\mu_1})$ for any $u_0 \in H$. Hence, for $\mu_1$-a.e. $u_0\in H$, $\mathbb{P}_{u_0}^\mathbb{N}(A_{\mu_1,u_0}) = \mathbb{P}_{u_0}^\mathbb{N}(A_{\mu_1}) = 1$. Similarly, we can find that $\mathbb{P}_{v_0}^\mathbb{N}(A_{\mu_2, v_0}) = 1$ for $\mu_2$-a.e. $v_0\in H$. Now, we choose a pair $u_0,v_0\in H$ such that $\mathbb{P}_{u_0}^\mathbb{N}(A_{\mu_1,u_0}) = \mathbb{P}_{v_0}^\mathbb{N}(A_{\mu_2,v_0}) =1$. Since $\mathbb{P}_{v_0} \sim \tilde{\mathbb{P}}_{v_0}$, by\;\Cref{lemma:girsanov}, we have that $\mathbb{P}_{v_0}^\mathbb{N}\sim \tilde{\mathbb{P}}_{v_0}^\mathbb{N}$. It follows that $\tilde{\mathbb{P}}_{v_0}(A_{\mu_2, v_0}) = 1$ for $\mu_2$-a.e. $v_0 \in H$. 

    Now, we let $\tilde{\xi}$ denote the law of ${(u(n), \tilde{v}(n))}_{n\in \mathbb{N}}$ where $u$ solves\;\eqref{eq:af} and $\tilde{v}$ solve\;\eqref{eq:nudgedaf} with initial conditions $u_0, v_0$ respectively. Define the set $D$ by 
    \[D := \{(u,\tilde{v}) \in H^\mathbb{N} \times H^\mathbb{N} : \lim_{n\rightarrow \infty} |u(n) - \tilde{v}(n)| = 0\}.\]
    By\;\Cref{corollary:posprob}, we have that
    \begin{equation}\label{eq:xi}
      \tilde{\xi}(D) = \mathbb{P}(\lim_{n\rightarrow \infty} |u(n) - \tilde{v}(n)| = 0) > 0.  
    \end{equation} 
    Now, define $\overline{D} := D \cap (A_{\mu_1 , u_0} \times A_{\mu_2, v_0})$. By\;\eqref{eq:xi}, we have $\tilde{\xi}(\overline{D}) > 0$ and hence nonempty. Choosing any $(u,\tilde{v}) \in \overline{D}$, we observe  
    \begin{align*}
        \bigg|\int_H \phi\ d\mu_1 - \int_H \phi\ d\mu_2\bigg| &= \bigg|\lim_{n\rightarrow \infty} \frac{1}{n} \sum_{k=1}^n \phi(u(k)) - \phi(\tilde{v}(k))\bigg|\\
        &\leq \limsup_{n\rightarrow \infty} \bigg|\frac{1}{n}\sum_{k=1}^n \phi(u(k)) - \phi(\tilde{v}(k))\bigg|\\
        &\leq \ell_{\phi} \limsup_{n\rightarrow \infty} \frac{1}{n}\sum_{k=1}^n |u(k) - \tilde{v}(k)|\\
        & = 0.
    \end{align*}
    Here $\ell_\phi$ is the Lipschitz constant for $\phi$. Since the left side of the inequality is lower bounded by $0$, we conclude that 
    \[\int_H \phi\ d\mu_1 = \int_H \phi\ d\mu_2.\]
    Since $\phi$ was arbitrary in $\mathcal{G}$ and $\mathcal{G}$ is a determining set, we have that $\mu_1 = \mu_2$ as desired. Combining the existence that was establish in\;\Cref{thm:erginv} with this result, we have uniqueness of an ergodic invariant measure.
\end{proof}

\section{Extensions}

\subsection{Convolution-Regularized 3D NSNE}
Consider the three-dimensional nonlinear stochastic Navier-Stokes equation on a bounded domain with a convolution-regularized inertial term of the form $B(k * u, u)$. That is, the equation
\[u(t) + \int_0^t [\nu_0Au + \nu_1\mathcal{A}_1(u) + B(k*u,u)]ds = u_0 + \int_0^t f ds + \int_0^t \sigma(u) dW(s)\]
where
\[k*u(x) = \int_{\mathbb{R}^3} k(x-y)u(y) dy,\quad k\in L^{6/5}(G).\]
Note that if the divergence-free and boundary conditions are in place, property\;\eqref{eq:tl} still holds. As a result, $\langle B(k*u,v),v \rangle = 0$ holds. Thus, a priori estimates similar to\;\Cref{prop:apriori} as well as moment estimates similar to\;\Cref{prop:moment} and\;\Cref{prop:momentnudged} follow. Since $k\in L^{6/5}(G)$, then, by the H\"{o}lder and the Young inequality for convolutions, we have
\begin{equation}\label{eq:regularized}
\begin{aligned}[b]
    |\langle B(k*u,u) - B(k*v,v), w \rangle| &= |\langle B(k*w,u),w\rangle|\\
    &\leq ||k*w||_{L^6} ||u||\, ||w||_{L^3}\\
    &\leq ||k||_{L^{6/5}} ||w||_{L^3}^2||u||\\
    &\leq ||k||_{L^{6/5}} ||w||\, |w|\, ||u||
\end{aligned}
\end{equation}
where $w = u-v$ and the last inequality is due to the fact that, in $\mathbb{R}^3$, $|w|_{L^3} \leq |w|^{1/2} ||w||^{1/2}$. Using estimate\;\eqref{eq:regularized} we arrive at estimates similar to those in\;\eqref{eq:gronwall} and\;\eqref{eq:diffN}. Hence, the results of this work apply to the convolution-regularized three-dimensional nonlinear stochastic Navier-Stokes equations.

\subsection{Unbounded Poincar\'{e} Domain}

Suppose $G$ is an unbounded Poincar\'{e} domain. That is, there exists a constant $\Lambda_P > 0$ such that, for all $\varphi \in W_0^{1,2}(G)$,
\begin{equation}\label{poincare}
    ||\varphi||_{L^2}^2 \leq \Lambda_P ||\nabla \varphi||_{L^2}^2.
\end{equation}
Since $G$ is a Poincar\'{e} domain, the embeddings\;\eqref{eq:embedding} remain to be compact, dense, and continuous. Note that the lower bound in\;\Cref{prop:lowerbound} does not hold since $G$ is unbounded. So, the estimates relying on\;\Cref{prop:lowerbound} do not follow automatically. However, since $|u|^{p-2}||u||^q$ is nonnegative and appears on the left side of\;\eqref{eq:eq1},\;\eqref{eq:eqp},\;\eqref{eq:nudgedeq1}, and\;\eqref{eq:nudgedpdif}, it can simply be removed to help obtain the moment estimates of the solution $u$ of\;\eqref{eq:af} and the solution $v$ of\;\eqref{eq:nudgedaf}. As a result, the estimates and constants no longer depend on $\nu_1$. Thus, the final results follow as before. 

\subsection{Other Viscosities}
Our main results rely on the energy equality and the associated moment estimates obtained for the solutions to the original equation and the nudged equation. In particular, these estimates provide the necessary bounds on the solution required for the long-time behavior and ergodic properties. Therefore, to extend the results to other types of viscosities, it is sufficient to verify that the fundamental energy structure is preserved and to provide the essential moment estimates on the solutions.\\

In this section, we consider another type of viscosity,
shear-dependent viscosity. The arguments required for this case differ only in the estimates involving the viscosity term, while the remaining components of the analysis follow as in the previous sections. For completeness, we provide the necessary modifications to the functional setting and establish the moment estimates of the solution in the subsection that follows.
\subsubsection{Shear-Dependent Viscosity}
Consider the stochastic Navier-Stokes equations with shear-dependent viscosity;
\begin{align}
    \frac{\partial u}{\partial t} - \nu_0 \Delta u -\nu_1||\nabla u||_{L^2(G)}^2 \Delta u + (u \cdot \nabla) u + \nabla \rho &= f(t) + \sigma(u) \frac{dW(t)}{dt} \notag\\
    \nabla\cdot u &= 0 \label{eq:sheardep}
\end{align}
with $u(x,t) = 0$ for $x\in \partial G$ and $u(x,0) = u_0(x)$ for all $x\in G$. The functional setting of this equation is similar to that of the usual stochastic Navier-Stokes equations. That is,
\begin{align*}
    H &= \{u\in L^2(G); \nabla \cdot u = 0, u \cdot \textbf{n}|_{\partial G} = 0\}\\
    V &= \{u \in W^{1,2}_0(G); \nabla \cdot u =0\}
\end{align*}
where $\mathbf{n}$ is the outward normal and $u\cdot \mathbf{n}\big|_{\partial G}$ is in the sense of trace operators. Let $H'$ and $V'$ denote the duals of $H$ and $V$ respectively. The norms of these spaces satisfy the following equivalences
\begin{align*}
    |u|^2 &:= ||u||_H^2 = \int_G |u|^2\, dx\\
    ||u||^2 &:= ||u||_V^2 = \int_G |\nabla u|^2\, dx
\end{align*}
The operators $A: V \rightarrow V'$ and $B: V \times V \rightarrow V'$ and the noise are as they are defined in Section 2. In addition to these, we define the nonlinear viscosity operator $\mathcal{A}_2(u)$ such that, for each $u \in V$,
\begin{align*}
    \langle \mathcal{A}_2(u),v\rangle = ||u||^2\int_G \nabla u \, \nabla v \, dx \quad \forall v\in V.
\end{align*}
It is clear that $\langle \mathcal{A}_2(u),v\rangle$ is linear in $v$. By the H\"{o}lder inequality, 
\begin{align*}
    |\langle \mathcal{A}_2(u),v\rangle| &= ||u||^2\Big|\int_G \nabla u\, \nabla v\, dx\Big|\\
    &\leq ||u||^3\, ||v|| < \infty,
\end{align*}
so the integral is well-defined. It follows that $\mathcal{A}_2(u)$ is bounded with $||\mathcal{A}_2(u)||_{V'} \leq ||u||^3$. Thus, the existence of such an operator $\mathcal{A}_2$ is ensured with $\mathcal{A}_2: V \rightarrow V'$ defined by $\mathcal{A}_2(u) = -\mathbf{P}_H ||u||^2 \Delta u$ where the equality holds in $V'$. This nonlinear viscosity operator satisfies the following property.
\begin{proposition}[\cite{s}, Theorem 1]
    The nonlinear viscous operator $\mathcal{A}_2$ is strictly monotone, i.e.
    \begin{equation}\label{eq:a2nl}
        \langle \mathcal{A}_2(u) - \mathcal{A}_2(v), u-v \rangle = \frac{1}{2}(||u||^2 + ||v||^2)||u-v||^2 + \frac{1}{2}{(||u||^2 - ||v||^2)}^2.
    \end{equation}
\end{proposition}
\begin{proof}
    Consider the integral
    \[\langle \mathcal{A}_2(u) - \mathcal{A}_2(v), u-v \rangle = \int_G (||u||^2 \nabla u - ||v||^2 \nabla v)\cdot (\nabla u - \nabla v) dx.\]
    This integral may be written as 
    \begin{align*}
        \frac{1}{2}(||u||^2 + ||v||^2) \int_G {(\nabla u - \nabla v)}^2 dx + \frac{1}{2}{(||u||^2 - ||v||^2)}^2 \int_G (\nabla u - \nabla v)\cdot (\nabla u + \nabla v) dx
    \end{align*}
    from which\;\eqref{eq:a2nl} follows.
\end{proof}
\noindent
Applying the Leray projection, we write\;\eqref{eq:sheardep} as 
\begin{equation}\label{eq:shearaf}
    u(t) + \int_0^t (\nu_0 Au + \nu_1\mathcal{A}_2(u) + B(u))ds = u_0 + \int_0^t f\, ds + \int_0^t \sigma(u) dW(s).
\end{equation}
Given the estimates on $A$, $B$, and the noise provided in Section 2 and the estimate on $\mathcal{A}_2$ mentioned here, the nonlinear operator $F$ defined by $F(u) = -\nu_0 Au - \nu_1\mathcal{A}_2(u) - B(u)$ satisfies a local monotonicity property and difference estimate similar to that of\;\Cref{prop:localmono} and\;\eqref{eq:gronwall} respectively. Thus, existence and uniqueness of a strong solution follow through similar techniques as those presented in\;\Cref{thm:solution}. The following moment estimates along the different estimate ensure the existence of an ergodic invariant measure.
\begin{proposition}\label{prop:shearmoments}
    Let $\sigma$ satisfy\;\Cref{ass:cv} and assume 
    \[\frac{\nu_0\lambda_1}{2} + 2\nu_1\lambda_1^2 > K.\]
    For $p\in \Big[2, 1 + \frac{\nu_0\lambda_1 + 4\nu_1\lambda_1^2}{2K}\Big)$, the solution $u$ of\;\eqref{eq:sheardep} satisfies
    \[\mathbb{E}[{|u(t)|}^p] \leq |u_0|^{p}e^{-at} + \frac{C}{a}\quad \forall t>0\]
    where
    \[a =\frac{p\nu_0\lambda_1}{4} + p\nu_1\lambda_1^2 - \frac{Kp(p-1)}{2}\]
    and 
    \[C = \begin{cases}
        \frac{2}{3\nu_0}||f||_{V'}^2 + K + 2\nu_1\lambda_1^2 & \text{if } p = 2\\
        \frac{2^{\frac{p}{2}-1}}{{\nu_0}^{p-1}{\lambda_1}^{\frac{p}{2}-1}}||f||^p_{V'} + \frac{{(p-1)}^{\frac{p}{2}} {(p-2)}^{\frac{p}{2}-1}}{{(\nu_0\lambda_1)}^{\frac{p}{2}-1}}K^{\frac{p}{2}} + p\nu_1\lambda_1^2 \quad & \text{if } p \in (2, 1+ \frac{\nu_0\lambda_1 + 4\tilde{C}\nu_1}{2K})
    \end{cases}.\]
\end{proposition}
\begin{proof}
   Applying the It\^{o} formula for $|u|^2$, we arrive at 
   \begin{align*}
    |u(t)|^2 + 2\nu_0\int_0^t ||u||^2 ds &+ 2\nu_1\int_0^t ||u||^4 ds\\
    &= |u_0|^2 + 2\int_0^t \langle f, u\rangle ds + \int_0^t ||\sigma(u)||_{L_Q}^2 ds + 2\int_0^t \langle u, \sigma(u)\; dW(s)\rangle.
   \end{align*}
   Observe that 
   \[2\nu_1||u||^4 \geq 2\nu_1\lambda_1^2(|u|^2 - 1).\]
   Applying the Young inequality, the assumption on $\sigma$, and the estimate above, we get
   \begin{align*}
        |u|^2 + \bigg(\frac{\nu_0\lambda_1}{2} + 2\nu_1\lambda_1^2 - K\bigg)\int_0^t |u|^2 ds \leq |u_0|^2 + \bigg(\frac{2}{3\nu_0}||f||_{V'}^2 + K + 2\nu_1\lambda_1^2\bigg)t + 2\int_0^t \langle u, \sigma(u) dW(s)\rangle.
   \end{align*}
   Taking the expectation and applying the Gr\"{o}nwall inequality, we arrive at 
   \[\mathbb{E}|u|^2 \leq |u_0|^2 e^{-at} + \frac{C}{a}\]
   where 
   \begin{align*}
        a&= \frac{\nu_0\lambda_1}{2} + 2\nu_1\lambda_1^2 - K\\
        C&= \frac{2}{3\nu_0}||f||_{V'}^2 + K + 2\nu_1\lambda_1^2.
   \end{align*}
   The proof for $p \in \Big(2, 1 + \frac{\nu_0\lambda_1 + 4\nu_1\lambda_1^2}{2K}\Big)$ follows similar to that of\;\Cref{prop:moment}. Combining the estimates on terms in\;\Cref{prop:moment} with the estimate $p\nu_1|u|^{p-2}||u||^4 \geq p\nu_1\lambda_1^2(|u|^p - 1)$ for the nonlinear viscosity term, we arrive at the desired moment estimates for $p \in \Big(2, 1 + \frac{\nu_0\lambda_1 + 4\nu_1\lambda_1^2}{2K}\Big)$.
\end{proof}
For uniqueness of the ergodic invariant measure, we define the nudged and stopped nudged equations in the same way as equations\;\eqref{eq:nudgedaf} and\;\eqref{eq:stoppednudged}. The existence of the generalized coupling is dependent only on the assumptions related to the stopped nudged term, so such a coupling still exists in this context. The only result that is needed to confirm the uniqueness of an ergodic invariant measure is the estimate in probability similar to\;\Cref{prop:estinprob} and the moment estimates of the nudged equation. From the energy equality, 
\begin{align*}
    |u(t)|^2 + 2\nu_0\int_0^t ||u||^2 ds &+ 2\nu_1\int_0^t ||u||^4 ds\\
    &= |u_0|^2 + 2\int_0^t \langle f, u\rangle ds + \int_0^t ||\sigma(u)||_{L_Q}^2 ds + 2\int_0^t \langle u, \sigma(u)\; dW(s)\rangle
\end{align*}
and the moment estimates of the solution to\;\eqref{eq:shearaf}, we can establish an estimate in probability similar to that of\;\Cref{prop:estinprob}. It remains to show moment estimates of the nudged equation. However, it is clear that the moment estimates of the nudged equation are derived using similar techniques to the moment estimates of the original equation. So, uniqueness of the ergodic invariant measure is guaranteed in this context.\\

\subsection{Other Models}
The results of this work are not specific to the stochastic Navier-Stokes equations, but rather follow from the underlying variational framework and the structural properties of the operators involved. In particular, the analysis relies on the coercivity of the dissipative operator, the cancellation properties of the nonlinear term, suitable continuity or local monotonicity estimates, and appropriate assumptions on the stochastic forcing. These techniques therefore extend to other nonlinear stochastic hydrodynamical models that possess the same analytical structure. We illustrate this by considering two examples: the stochastic shell model and the stochastic magnetohydrodynamics (MHD) system.
\subsubsection{Shell Model of Turbulence}
For $n= 1, 2,\dots$ and $u_n \in \mathbb{C}$, consider the system
\begin{equation}\label{eq:shell}
    \frac{du_n}{dt} + \nu k_n^2 u_n + i(ak_{n+1}u_{n+1}^*u_{n+2}^* + bk_n u_{n-1}^*u_{n+1}^* + ck_{n-1}u_{n-1}^*u_{n-2}^*) = f_n
\end{equation}
with boundary conditions $u_{-1}= u_0 = 0$. Here $\nu$ denotes the kinematic viscosity, $u_n^*$ the complex conjugate of $u_n$, and $a,b,c$ are real-valued parameters such that the conservation of energy is preserved: $a+b+c=0$. The standard model takes $a=-1$, $b = \frac{1}{2}$, and $c = \frac{1}{2}$. The shell model that is described by equation\;\eqref{eq:shell} is due to Gledger, Ohkitani, and Yamada\;\cite{goy} and is known as the GOY model.\\

This model arises when the Navier-Stokes equation is written in the Fourier space where the interaction between modes is preserved between the nearest modes. For the GOY model, it describes a one-dimensional cascade of energies among an infinite sequence of complex velocities, $u_n(t)$, on a one-dimensional sequence of wave numbers, $k_n = k_0 2^n$, for $n = 1, 2, \dots$ where $k_0 > 0$ and the index $n$ is referred to as the ``shell index'' (see\;\cite{turb}). The wave numbers, $k_n$, may also be taken as $k_n = k_0 \lambda^n$ for $\lambda > 1$.\\

Now, we establish the functional setting of\;\eqref{eq:shell}. Define the real Hilbert space
\[H := \Big\{u = (u_1, u_2,\dots) \in \mathbb{C}^\infty : \sum_{n=1}^\infty |u_n|^2 < \infty\Big\}\]
endowed with scalar product $\langle \cdot, \cdot \rangle_H$ and norm $|\cdot |$ given by 
\[\langle u, v\rangle_H = Re \sum_{n=1}^\infty u_n v_n^*, \qquad |u|:= ||u||_H = {\bigg(\sum_{n=1}^\infty |u_n|^2\bigg)}^{1/2}\]
for every $u,v \in H$ and define the space 
\[V := \Big\{u \in H : \sum_{n=1}^\infty k_n^2 |u_n|^2 < \infty\Big\}\]
with norm 
\[||u|| := ||u||_V = {\bigg(\sum_{n=1}^\infty k_n^2 |u_n|^2\bigg)}^{1/2}.\]
The operator $A: D(A) \rightarrow H$ defined by 
\[Au = (k_1^2 u_1, k_2^2 u_2, \dots ), \quad \forall u \in D(A)\]
is a linear operator that is positive definite and self-adjoint. Since $A$ is positive definite, we are able to define the power $A^{1/2}$ where 
\[A^{1/2}u = (k_1 u_1, k_2 u_2, \dots), \quad \forall u = (u_1, u_2, \dots).\]
Note that, when equipped with the scalar product and norm
\[\langle u,v\rangle_{D(A^{1/2})} = \langle A^{1/2} u, A^{1/2} v\rangle, \quad ||u||_{D(A^{1/2})} = {\bigg(\sum_{n=1}^\infty k_n^2 |u_n|^2\bigg)}^{1/2} \quad \forall u, v \in D(A^{1/2}),\]
we have $V = D(A^{1/2})$. If we consider the dual of $V$ as $V' = D(A^{-1/2})$, then the embeddings 
\[V\hookrightarrow H = H' \hookrightarrow V'\]
is dense, continuous, and compact. We now define the bilinear operator $B: V \times H \rightarrow H$ as 
\[B(u,v) = (B_1(u,v), B_2(u,v), \dots)\]
where 
\[B_n(u,v) = ik_n\Big(\frac{1}{4} u_{n+1}^* v_{n-1}^* - \frac{1}{2}(u_{n+1}^* v_{n+2}^* + u_{n+2}^* v_{n+1}^*) + \frac{1}{8}u_{n-1}^* v_{n-2}^*\Big)\] 
Let ${\{e_n\}}_{n=1}^\infty$ denote the standard orthonormal basis of $H$. Then the bilinear operator $B$ can be written as 
\[B(u,v) = i\sum_{n=1}^\infty k_n\Big(\frac{1}{4} u_{n+1}^* v_{n-1}^* - \frac{1}{2}(u_{n+1}^* v_{n+2}^* + u_{n+2}^* v_{n+1}^*) + \frac{1}{8}u_{n-1}^* v_{n-2}^*\Big)e_n.\]
The following lemma provides some estimates and the cancellation property of $B$ (see\;\cite{bbbf}).
\begin{lemma}
    There is a constant $C>0$ such that
    \begin{equation}
        |B(u,v)| \leq C ||u||\, |v|, \quad \forall u\in V, v\in H
    \end{equation}
    and 
    \begin{equation}
        |B(u,v)|\leq C |u|\, ||v||, \quad \forall u\in H, v\in V.
    \end{equation}
    Also, for every $u\in V$ and $v\in H$, we have 
    \begin{equation}\label{eq:goycancel}
        \langle B(u,v),v\rangle_H = 0.
    \end{equation}
\end{lemma}
With the functional setting well established, we introduce the stochastic shell equation for the GOY model on $H$
\begin{equation}
    du + [\nu Au + B(u,u)]dt = f(t) dt + \sigma(u) dW(t)
\end{equation}
where $u= (u_1,u_2,\dots)$ and the operators $A$ and $B$ are defined above. Utilizing the bilinearity of $B$ and the cancellation property\;\eqref{eq:goycancel}, we say that $\langle B(u,u) - B(v,v), w\rangle = \langle B(w,v), w\rangle$ where $w = u-v$. It follows, by this fact and the well-known estimate $||w||_{\ell^4}\leq |w|^{1/2} ||w||^{1/2}$, that 
\begin{align*}
    |\langle B(u,u) - B(v,v), w\rangle| &= |\langle B(w,v), w\rangle|\\
    &\leq ||w||_{\ell^4} ||w||\, ||v||_{\ell^4}\\
    &\leq |w|^{1/2} ||w||^{3/2} ||v||_{\ell^4}\\
    &\leq \frac{1}{2} ||w||^2 + \frac{27}{32} |w|^2 ||v||_{\ell^4}^4.
\end{align*}
From here, the proof of local monotonicity follows from similar steps as those in\;\Cref{prop:localmono}, so existence and uniqueness of a strong solution is assured. Also, due to the cancellation property of $B$ and the coercivity of $A$, we are able to derive essential estimate in probability and the moment estimates need to prove existence and uniqueness of an ergodic invariant measure.

\subsubsection{Magnetohydrodynamical (MHD) System}
Let $G\subset \mathbb{R}^2$ be a bounded, open, and simply connected domain with smooth boundary $\partial G$. The stochastic MHD system (see\;\cite{mhd}) in the non-dimensional form is given by
\begin{equation}\label{eq:mhd}
    \frac{\partial u}{\partial t} + (u\cdot \nabla) u - \nu\Delta u- S(B\cdot \nabla)B + \nabla\Big(p + \frac{S|B|^2}{2}\Big) = 0.
\end{equation}
Here $u$ is the velocity, $B$ the magnetic field, $p$ the pressure field, and $\nu$ the viscosity. Note that the viscosity is inversely related to the Reynolds number $R_e$, i.e. $\nu = \frac{1}{R_e}$. Let $S := \frac{M^2}{R_e R_m}$ where $M$ is the Hartman number and $R_m$ is the magnetic Reynolds number. The Maxwell equation is
\begin{equation}\label{eq:maxwell}
    \frac{\partial B}{\partial t} + (u \cdot \nabla)B + \eta \text{ curl}(\text{curl } B) - (B \cdot \nabla) u = 0
\end{equation}
where $\nabla \cdot u = 0$ and $\nabla \cdot B = 0$. Here $\eta = \frac{1}{R_m}$. Along with equations\;\eqref{eq:mhd} and\;\eqref{eq:maxwell}, we have the boundary conditions $u = 0$, $B \cdot \textbf{n} = 0$, and $\text{curl } B = 0$ on $\partial G$ where $\textbf{n}$ is the unit outer normal on $\partial G$ and initial conditions $u(x,0) = u_0(x)$ and $B(x,0) = B_0(x)$ for $x \in G$. Note that if $u = (u_1, u_2)$ is a vector-valued function defined on $G$, then $\text{curl }u = \frac{\partial u_2}{\partial x} - \frac{\partial u_1}{\partial y}$. Also,  if $\varphi$ is a scalar-valued function defined on $G$, then $\text{curl } \varphi = (\frac{\partial \varphi}{\partial y}, -\frac{\partial \varphi}{\partial x})$.\\

For the functional setting, we set $H = H_1 \times H_2$ and $V = V_1 \times V_2$ where 
\begin{align*}
    H_1 &= \{\varphi \in L^2(G) : \nabla \cdot \varphi = 0,\; \varphi \cdot \textbf{n}|_{\partial G} = 0\}\\
    H_2 &= H_1\\
    V_1 &= \{\varphi \in W^{1,2}_0(G) : \textbf{n} \cdot \varphi = 0\}\\
    V_2 &= \{\varphi \in W^{1,2}(G) : \textbf{n} \cdot \varphi = 0,\; \varphi \cdot \textbf{n} |_{\partial G} = 0\}.
\end{align*}
The spaces $H_1$ and $H_2$ are equipped with the $L^2(G)$ norm. Define the inner product by $[y_1,y_2] = \langle u_1, u_2\rangle_{H_1} + S\langle B_1, B_2\rangle_{H_2}$ where $y_i = \{u_i, B_i\}$. This inner product is equivalent to $\langle y_1, y_2\rangle_H = \langle u_1, u_2\rangle_{H_1} + \langle B_1, B_2\rangle_{H_2}$. The norm on $H$ is given by $||y||_H = \sqrt{\langle y, y\rangle_H}$. We endow $V_1$ with the inner product $\langle \varphi, \psi \rangle_{V_1} = \langle \nabla \varphi, \nabla \psi \rangle_{L^2(G)}$ and norm $||\varphi||_{V_1} = \sqrt{\langle \varphi, \varphi\rangle}$. Note that the $V_1$ norm is equivalent to the $W^{1,2}(G)$ norm due to the Poincar\'{e} inequality. The space $V_2$ is endowed with the inner product $\langle \varphi, \psi\rangle_{V_2} = \langle \text{curl }\varphi, \text{ curl } \psi\rangle_{L^2(G)}$ and the norm $||\varphi||_{V_2} = \sqrt{\langle \text{curl }\varphi, \text{ curl } \varphi\rangle_{L^2(G)}}$. The norm of $V_2$ is also equivalent to the $W^{1,2}(G)$ norm (see\;\cite{tm}). We now assign the space $V$ the scalar product given by $[|y_1,y_2|] = \langle u_1, u_2 \rangle_{V_1} + S\langle B_1, B_2\rangle_{V_2}$. With the function spaces well-established, we have the embedding
\[V \hookrightarrow H = H' \hookrightarrow V'\]
which is dense, continuous, and compact.\\

The operator $A: V\rightarrow V'$ is defined by 
\[\langle Ay_1, y_2\rangle = a(y_1,y_2)\quad \text{with } a(y_1,y_2) = \nu\langle u_1, u_2\rangle_{V_1} + S\eta\langle B_1, B_2 \rangle_{V_2} \quad \forall y_i = \{u_i, B_i\} \in V.\]
Hence,
\[\langle Ay, y\rangle = \nu||\nabla u||_{L^2(G)}^2 + S\eta||\text{curl } B||_{L^2{G}}^2 \geq c||y||_H^2\]
for some positive constant $c$ depending on $\nu, \eta$, and $S$. Thus, the linear operator $A$ is symmetric, coercive, and positive definite.\\

Define the operator $B: V \times V \rightarrow V'$ as the continuous bilinear operator such that 
\[\langle B(y_1, y_2), y_3\rangle = b(y_1, y_2, y_3)\]
where the trilinear operator $b: V\times V\times V \rightarrow \mathbb{R}$ is given by 
\[b(y_1, y_2, y_3) = \tilde{b}(u_1, u_2, u_3) - S\tilde{b}(B_1, B_2, u_3) + S\tilde{b}(u_1, B_2, B_3) - S\tilde{b}(B_1, u_2, B_3) \]
for all $y_i = \{u_i, B_i\} \in V$ and $\tilde{b}: {(W^{1,2}(G))}^{\otimes 3} \rightarrow \mathbb{R}$ is defined by 
\[\tilde{b}(\varphi, \psi, \theta) = \sum_{i,j = 1}^2 \int_{G} \varphi_i \frac{\partial \psi_j}{\partial x_i} \theta_j\, dx.\]
By integration by parts, we have that 
\[\tilde{b}(\varphi, \psi, \psi) = \sum_{i,j = 1}^2 \int_{G} \varphi_i \frac{\partial \psi_j}{\partial x_i} \psi_j\, dx = 0\]
for all $\varphi \in W^{1,2}(G)$ and $\psi \in V_1$. Additionally, we have $\tilde{b}(\varphi, \psi, \theta) = -\tilde{b}(\varphi, \theta, \psi)$ if $\psi$ or $\theta \in V_1$ and $\varphi \in W^{1,2}(G)$. It follows that $ \langle B(y,z), z\rangle = b(y,z,z) = 0$.\\

With operators $A$ and $B$ now well defined, the stochastic MHD system appears as 
\begin{equation}\label{eq:mhdaf}
    dy + [Ay + B(y,y)]dt = f(t) dt + \sigma(y)dW(t)
\end{equation}
with $y(0)\in H$. To show that\;\eqref{eq:mhdaf} satisfies the local monotonicity property, we need to show an estimate on the bilinear term $B$ similar to that of\;\eqref{eq:bw}. We show such an inequality now. Let $\hat{y} = y_1 - y_2$, $\hat{u} = u_1 - u_2$, and $\hat{B} = B_1 - B_2$. For simplicity of notation below, we let $|\cdot |$ and $||\cdot||$ denote the appropriate $H$ and $V$ respectively depending on context. For example, $|y_1|$ denotes the $H$ norm while $|u_1|$ would denote the $H_1$ norm. By definition of $B$ and various known inequalities, we get
\begin{align*}
|\langle B(y_1)-B(y_2),\hat{y}\rangle|&=|b(\hat{y},y_1,\hat{y})|\\
&=\Big|\tilde{b}(\hat{u},u_1,\hat{u})-S\tilde{b}(\hat{B},B_1,\hat{u})+S\tilde{b}(\hat{u},B_1,\hat{B})-S\tilde{b}(\hat{B},u_1,\hat{B})\Big|\\
&\leq|\tilde{b}(\hat{u},u_1,\hat{u})|+S|\tilde{b}(\hat{B},B_1,\hat{u})|+S|\tilde{b}(\hat{u},B_1,\hat{B})|+S|\tilde{b}(\hat{B},u_1,\hat{B})|\\
&\leq C|\hat{u}|^{1/2}\|\hat{u}\|^{1/2}\|u_1\||\hat{u}|^{1/2}\|\hat{u}\|^{1/2} +CS\|B_1\||\hat{B}|^{1/2}\|\hat{B}\|^{1/2}|\hat{u}|^{1/2}\|\hat{u}\|^{1/2}\\
&\quad +CS\|B_1\||\hat{u}|^{1/2}\|\hat{u}\|^{1/2}|\hat{B}|^{1/2}\|\hat{B}\|^{1/2} +CS\|u_1\||\hat{B}|^{1/2}\|\hat{B}\|^{1/2}|\hat{B}|^{1/2}\|\hat{B}\|^{1/2}\\
&\leq C\|u_1\||\hat{u}|\|\hat{u}\| +CS\|B_1\|{(|\hat{u}|\|\hat{u}\||\hat{B}|\|\hat{B}\|)}^{1/2} +CS\|u_1\| |\hat{B}|\|\hat{B}\|\\
&\leq C\|u_1\||\hat{u}|\|\hat{u}\|+CS\|B_1\|(\varepsilon |\hat{u}|\|\hat{u}\|+\frac{1}{4\varepsilon}|\hat{B}|\|\hat{B}\|)+CS\|u_1\| |\hat{B}|\|\hat{B}\|\\
&\leq \varepsilon(\|\hat{u}\|^2+\|\hat{B}\|^2)+\frac{C}{\varepsilon}(\|u_1\|^2+\|B_1\|^2)(|\hat{u}|^2+|\hat{B}|^2)\\
&= \varepsilon ||\hat{y}||^2 + \frac{C}{\varepsilon}||y_1||^2 |\hat{y}|^2.
\end{align*}
From here, the local monotonicity property follows in a straightforward way, so existence and uniqueness of a strong solution is assured. The energy equality is as follows:
\[d|y|^2 + 2\langle Ay,y\rangle dt = 2\langle f,y\rangle dt + ||\sigma(y)||_{L(Q)}^2 + 2\langle y, \sigma(y) dW(t)\rangle.\]
Through some work, similar to that of\;\Cref{prop:moment}, the moment estimates of the solution to\;\eqref{eq:mhdaf} can be achieved with
\[\mathbb{E}|y|^p \leq |y(0)|^p e^{-at} + \frac{C}{a}\]
for $a = \frac{p\lambda_A}{2} - p(p-1)K$ and $C = \frac{p(p-1)}{2}K + {\Big(\frac{2(p-1)}{p\lambda_A}\Big)}^{p-1}||f||_{V'}^p$ where $\lambda_A = \lambda_1 \min\{\nu, \eta\}$ and provided $\lambda_A > 2(p-1)K$ and $p \in [2, 1 + \frac{\lambda_A}{2K})$. From here, the proof of existence and uniqueness of an ergodic invariant measure follows from the techniques in this paper.

\section*{Acknowledgements}
The author would like to thank Professor P. Sundar for his guidance and helpful discussions during this research. The author is especially grateful for his valuable feedback and encouragement during the development of this work.

\end{document}